\documentclass[11pt,letterpaper]{article}
\usepackage[letterpaper,textwidth=455bp,textheight=650bp,centering]{geometry}
\usepackage{amsmath,amssymb,amsthm,hyperref,fancyhdr,needspace,etoolbox}
\hypersetup{hidelinks,pdfauthor={Shaoqiang Deng and Junna Ni},pdftitle={The modular function of a group quasi-isometric to a horospherical product},pdfsubject={Geometric group theory; MSC 2020: 20F65, 20F69, 22D05, 53C30, 20E08},pdfkeywords={quasi-isometric rigidity, horospherical product, modular function, locally compact group, Heintze space, amenable hyperbolic group, Diestel-Leader graph}}
\newcommand{\hor}{\mathbin{\triangleright\!\!\triangleleft}}
\newcommand{\R}{\mathbb R}\newcommand{\Z}{\mathbb Z}\newcommand{\Q}{\mathbb Q}
\newcommand{\tr}{\operatorname{tr}}\newcommand{\Isom}{\operatorname{Isom}}
\newcommand{\Aut}{\operatorname{Aut}}\newcommand{\DL}{\operatorname{DL}}\newcommand{\HT}{\operatorname{HT}}\newcommand{\Sol}{\mathrm{Sol}}
\newcommand{\one}[1]{\mathbf 1_{#1}}
\newtheorem{thm}{Theorem}[section]
\newtheorem{lem}[thm]{Lemma}\newtheorem{prop}[thm]{Proposition}\newtheorem{cor}[thm]{Corollary}
\theoremstyle{definition}\newtheorem{defn}[thm]{Definition}\newtheorem{ex}[thm]{Example}\newtheorem{rem}[thm]{Remark}
\theoremstyle{plain}\newtheorem{mainresult}{Theorem}
\newtheorem{maincor}[mainresult]{Corollary}
\AtBeginEnvironment{defn}{\Needspace{5\baselineskip}}
\AtBeginEnvironment{lem}{\Needspace{5\baselineskip}}
\AtBeginEnvironment{thm}{\Needspace{5\baselineskip}}
\AtBeginEnvironment{mainresult}{\Needspace{6\baselineskip}}
\AtBeginEnvironment{maincor}{\Needspace{5\baselineskip}}
\title{The Modular Function of a Group Quasi-isometric to a Horospherical Product}
\author{Shaoqiang Deng$^{1}$ and Junna Ni$^{2}$\\[5pt]
  {\small $^{1}$School of Mathematical Sciences and LPMC, Nankai University}\\
  {\small Tianjin 300071, P. R. China \quad \href{mailto:dengsq@nankai.edu.cn}{\texttt{dengsq@nankai.edu.cn}}}\\[3pt]
  {\small $^{2}$School of Mathematical Sciences, South China Normal University}\\
  {\small Guangzhou 510631, P. R. China \quad \href{mailto:nijunna@126.com}{\texttt{nijunna@126.com}}}}
\date{}
\begin{document}\maketitle
\begin{abstract}
Every compactly generated locally compact group quasi-isometric to a horospherical product of unequal growth exponents is amenable and non-unimodular. We determine its modular function by the identity $\log\Delta_H=-(m_X-m_Y)\chi$, where $\chi$ is the stable height character and the factors are Heintze spaces, regular metric trees or millefeuille spaces. Modulo a compact normal subgroup, the group acts faithfully and geometrically on a product of the same factor types, with the original exponents and a rough-isometric comparison. This excludes finitely generated groups and proves the non-unimodular case of the Dymarz--Fisher--Xie conjecture for SOL-like groups with two nontrivial factors. For Carnot and Carnot-by-Carnot factors, the equal-trace case also allows us to remove the additional hypothesis in their quasi-isometric rigidity theorem. A boundary theorem supplies the amenability used in the proof: the closure of a uniform quasi-similarity group acting cocompactly on distinct boundary pairs is locally compact and focal hyperbolic. It also yields automatic amenability in several Tukia-type theorems.
\end{abstract}


\section{Introduction}
\subsection{The problem and the main theorem}
Dymarz, Fisher and Xie ask whether a SOL-like group can be quasi-isometric to a finitely generated group when none of its left-invariant Riemannian metrics has an isometry group admitting a uniform lattice \cite[Conjecture~1.6]{DFX25}. We prove the conjecture for every non-unimodular SOL-like group with two nontrivial factors. More generally, for a class of horospherical products that includes trees and millefeuille spaces, we determine the modular function of every compactly generated locally compact group in the quasi-isometry class. The formula gives an obstruction to discrete groups and describes the locally compact groups that can occur.

Write $h_X,h_Y$ for heights increasing towards chosen points at infinity and put
\[
 Z=X\hor Y=\{(x,y):h_X(x)+h_Y(y)=0\},\qquad h(x,y)=h_X(x).
\]
An ascent in $X$ forces a descent of the same length in $Y$. Two hyperbolic planes give the geometry of $\Sol$, and two regular trees give the Diestel--Leader examples. We allow Heintze spaces, regular metric trees and millefeuille spaces, as specified in Definition~\ref{def:model}. Their horospherical growth exponents are $m_X,m_Y$.

The exponents measure the mass gained on descending through horospheres. In the upper half-plane, $z\mapsto e^t z$ moves $i$ upwards by height $t$ and multiplies Euclidean boundary distances by $e^t$. In a tree with $p$ children and edge length $\ell$, descent through $k$ edges gives $p^k$ descendants, so the rate is $(\log p)/\ell$. A shift in the product expands one horizontal direction and contracts the other. This suggests that the difference of the two rates should determine Haar scaling. A quasi-isometry initially supplies only a quasi-action, so neither an exact shift nor a geometric action is given.

Fix a left Haar measure $\lambda_H$, with convention
\[
                 \lambda_H(Eg)=\Delta_H(g)\lambda_H(E).
\]
For example, $(x,t)(y,s)=(x+e^t y,t+s)$ has left Haar measure $e^{-t}\,dx\,dt$ and $\Delta(0,t)=e^{-t}$. This fixes the sign below. Choose a quasi-isometry $\varphi:H\to Z$, a coarse inverse $\bar\varphi$ and $o\in Z$. Set $q_g=\varphi L_g\bar\varphi$ and $b_g=h(q_g o)-h(o)$. Following powers removes the bounded ambiguity in $b_g$. A \emph{rough isometry} is a coarsely onto map preserving distances up to a bounded additive error.
\begin{mainresult}\label{thm:A}
Let $X,Y$ be model spaces with $m_X\ne m_Y$, and let $H$ be a compactly generated locally compact group quasi-isometric to $Z=X\hor Y$.
\begin{enumerate}
\item The group $H$ is amenable. The limits
\[
     \chi(g)=\lim_{n\to\infty}\frac{b_{g^n}}{n}
\]
define a nonzero continuous homomorphism $H\to\R$, and
\begin{equation}\label{eq:modular}
                \log\Delta_H(g)=-(m_X-m_Y)\chi(g)\qquad(g\in H).
\end{equation}
In particular $H$ is non-unimodular.
\item There are a compact normal subgroup $W\lhd H$ and model spaces $X',Y'$ of the respective types of $X,Y$ with $m_{X'}=m_X$ and $m_{Y'}=m_Y$. The group $H/W$ acts faithfully and geometrically by pairs of isometries on $Z'=X'\hor Y'$. Each $g\in H$ shifts the two heights by $\chi(g)$ and $-\chi(g)$, respectively. There are rough isometries $X\to X'$, $Y\to Y'$ and a rough isometry $F:Z\to Z'$ with
\[
                  \sup_{g\in H,\ p\in Z}d_{Z'}(Fq_g p,gFp)<\infty .
\]
\end{enumerate}
\end{mainresult}
Equation~\eqref{eq:modular} identifies the stable height directly from $H$ and the two exponents. In particular it is independent of the quasi-isometry, the coarse inverse and the base point, once the factors have been ordered and their height unit fixed. Reversing the order changes the signs of both $\chi$ and $m_X-m_Y$. A common rescaling changes the numerical exponents and the height by reciprocal factors.
\begin{rem}[Metric amenability and the quantitative conclusion]\label{rem:metricamenability}
For a compactly generated locally compact group, metric amenability is equivalent to amenability together with unimodularity, and is a quasi-isometry invariant \cite{Te08}, \cite[Lemma~19.5]{Cor18}. Thus, once amenability of $H$ is established, a separate proof that a group model for $Z$ is not metrically amenable would also imply non-unimodularity of $H$. This criterion alone does not establish amenability of $H$: ordinary amenability is not a quasi-isometry invariant among locally compact groups. Equation~\eqref{eq:modular} gives additional quantitative information, identifying the entire modular function and the normalized height character. Our proof obtains this information from the boundary actions and their geometric realizations; Appendix~\ref{sec:counting} gives a separate proof of the qualitative obstruction.
\end{rem}

\subsection{Relation to Ferragut's geometric rigidity theorem}
Ferragut proves that quasi-isometries between unequal-exponent horospherical products lie at uniformly bounded distance from product maps \cite[Theorem~4.1]{Fe26b}. His ratio invariant $m_X/m_Y$ is stated in \cite[Theorem~B]{Fe26b}. His work extends the coarse-differentiation methods of Eskin, Fisher and Whyte \cite{EFW12}. We use this factorization as the input for constructing two boundary actions of $H$.

The group-theoretic conclusions should be distinguished from this geometric input. In the smooth case, Ferragut classifies non-unimodular Carnot--Sol Lie groups up to quasi-isometry \cite[Theorem~C]{Fe26b} and describes the self-quasi-isometry group at unequal traces \cite[Theorem~E]{Fe26b}. Dymarz, Fisher and Xie extend the isomorphism classification within the class of SOL-like groups with Carnot or Carnot-by-Carnot factors \cite[Theorem~1.8]{DFX25}. These classification statements compare members of specified Lie-group classes. In Theorems~\ref{thm:A} and~\ref{thm:D}, the group $H$ is instead an arbitrary compactly generated locally compact group. The conclusions include a compact-kernel geometric realization and the modular function of $H$. In particular, the locally compact group assertion in Corollary~\ref{cor:C}(a) does not follow from Ferragut's Theorem~C alone. This distinction does not restrict Ferragut's factorization theorem to Lie groups: its admissible-space formulation is precisely the one used here.

Compact-kernel realization in an isometry group has precedents in \cite[Theorems~5.8(2) and~5.9(2)]{DX16}, for model filiform Heintze groups and certain abelian Heintze groups. We use the same general strategy of passing to boundary actions and conjugating them to similarities. The two-factor setting additionally requires compatibility of the two height characters; Theorem~\ref{thm:A} allows regular-tree and millefeuille factors as well as smooth factors.

To join geometric realizations of those actions, we need height and boundary control for the same factor quasi-isometry. Lemma~\ref{lem:height} records a boundary map $u$ and a height map $\phi$ for which
\[
 \log\rho_{X'}(u\xi,u\eta)=\phi(\log\rho_X(\xi,\eta))+O(1),\qquad
 \phi(s+a)-\phi(s)=\frac{m_X}{m_{X'}}a+O(1).
\]
The errors are uniform in the endpoints, the base height $s$ and the increment $a$. This formulation tracks the height unit through the boundary action for each of the three factor types. Lemma~\ref{lem:lift} also constructs rough-isometric lifts from individual boundary quasi-similarities. Section~\ref{sec:actions} compares the resulting coarse heights with exact Busemann characters, and Section~\ref{sec:mod} identifies these characters with the modular function. These steps supply the group-action and Haar-measure conclusions of Theorem~\ref{thm:A}.

\subsection{Consequences for SOL-like groups and graphs}
For $S=(N_+\times N_-)\rtimes\R$, with action $(e^{tD_+},e^{-tD_-})$, the exponents are proportional to $\tr D_+$ and $\tr D_-$. At unequal traces, Corollaries~\ref{cor:B} and~\ref{cor:C} show that every compactly generated locally compact group quasi-isometric to $S$ is amenable and non-unimodular. A discrete group's counting measure is right invariant, so this excludes all finitely generated groups. Applying the same conclusion to $\Isom(S,g)$ for each left-invariant Riemannian metric $g$ also excludes uniform lattices in these isometry groups. Thus both the premise and the conclusion of \cite[Conjecture~1.6]{DFX25} hold in this case.

For Carnot and Carnot-by-Carnot factors, Theorem~\ref{thm:D} realizes the geometric action on the original SOL-like group with a fixed left-invariant metric. In \cite[Theorem~1.5]{DFX25}, if a Carnot-by-Carnot factor has a one-dimensional ideal generated by its lowest eigenspace, the additional assumption is that the finitely generated group $\Gamma$ is amenable, or that some $\Isom(S,g)$ admits a uniform lattice. At equal traces, metric amenability and its quasi-isometry invariance imply amenability of $\Gamma$, verifying the first alternative. At unequal traces, Corollary~\ref{cor:B} excludes $\Gamma$ altogether. Corollary~\ref{cor:E} combines this obstruction with their theorem and the metric-amenability criterion; it closes the remaining case of their conjecture in this class of factors without requiring a further rigidity argument.

The discrete exclusion for classical $\Sol(m,n)$ and unequal-parameter Diestel--Leader graphs was proved by Eskin, Fisher and Whyte \cite[Theorems~1.2 and~1.4]{EFW12}. Pallier's appendix \cite[Theorem~A.1]{Pa23} already treats products of model filiform groups of nilpotency classes $3\le n<m$ with standard Carnot derivations. Their traces are $1+n(n+1)/2$ and $1+m(m+1)/2$, so this is a direct special case of Corollary~\ref{cor:C}(a). The Lie-group classification in \cite[Theorem~1.8]{DFX25} is another antecedent. The new assertions here concern the arbitrary locally compact group $H$, its modular function and its geometric realization, rather than these known discrete exclusions or classifications within specified Lie-group classes.

For arbitrary proper geodesic spaces in the unequal-exponent class, Corollary~\ref{cor:B} excludes lattices in every closed group acting cocompactly on the space. Corollary~\ref{cor:F} applies this to automorphism groups of quasi-transitive graphs in the unequal-parameter Diestel--Leader classes. Context for geometric realizations is provided by \cite{Dy15}, \cite{CFK12}, \cite{CL25} and \cite[Section~1.2.2]{LDPX25}.

\subsection{What happens at equal exponents}
At equality, the coefficient in~\eqref{eq:modular} vanishes. Lattices in $\Sol$, lamplighter groups and $\mathrm{BS}(1,p)$ show that discrete groups can occur. A mixed example of Caplinger and Levitin \cite[Section~1.2]{CL25} explains why the branching contribution to $m_X$ is essential. The matrix
\[
 A=\begin{pmatrix}3&1\\1&1\end{pmatrix}
\]
has eigenvalues $\alpha=2+\sqrt2$ and $\beta=2-\sqrt2$, with $\alpha\beta=2$. Its ascending HNN extension acts geometrically on a product of a Heintze plane and a mixed factor with a $3$-regular tree. In the height unit of Example~\ref{ex:models}(d),
\[
              m_X=\log\alpha,\qquad
              m_Y=-\log\beta+\log2=\log\alpha.
\]
The smooth rates differ by $\log2$, and the branching contributes exactly that amount. The full construction and height conventions are given in that example. Double-mixed examples also occur \cite[Section~5]{Cap26}. Thus the equal-exponent setting contains further questions about geometric actions with mixed factors, including the cases raised in \cite[Question~1]{CL25}.

\subsection{A boundary theorem}
A vertical line in a model space $X$ is specified by its lower endpoint $\zeta$ in the punctured boundary $\Omega_X$. We equip $\Omega_X$ with a quasi-distance $\rho_X$ such that $\log\rho_X(\zeta,\eta)$ records the height at which the two lines approach. A homeomorphism is an $L$-\emph{quasi-similarity} if it multiplies all distances by a number $r>0$ up to factors between $L^{-1}$ and $L$. In a uniform group of quasi-similarities, $r$ varies with the element while $L$ is fixed.

Think of two distinct points on the real line. Their location and their separation specify a position and a scale. In the hyperbolic model these determine a point up to bounded error. Cocompactness on distinct pairs therefore supplies the missing control on both horizontal displacement and vertical scale. A hyperbolic locally compact group with more than two boundary points and a global fixed point on its boundary is called \emph{focal}.
\setcounter{mainresult}{6}
\begin{mainresult}\label{thm:G}
Let $\Gamma$ be a group of $L$-quasi-similarities of $\Omega_X$ acting cocompactly on ordered pairs of distinct points. Its closure $Q$ in the homeomorphism group of the one-point compactification is locally compact, consists of $L$-quasi-similarities, and is a focal hyperbolic group quasi-isometric to $X$. In particular $Q$ is amenable.
\end{mainresult}
Amenability here is a property of the closure with its locally compact topology. For example, translations, dilations and a nonabelian free subgroup of $\mathrm{SO}(n)$, $n\ge3$, generate a group of similarities acting cocompactly on pairs whose underlying discrete group is non-amenable. Section~\ref{sec:tukia} applies Theorem~\ref{thm:G} to the hypotheses in \cite{DX16} and \cite{DFX25}.

\subsection{Two proofs of the non-unimodularity obstruction}
Ferragut's factorization gives the two genuine boundary actions in Sections~\ref{sec:models} and~\ref{sec:product}. Section~\ref{sec:closure} proves local compactness and focality of their closures. Section~\ref{sec:actions} realizes the two actions geometrically, compares their height units and proves that their exact shifts cancel. Section~\ref{sec:mod} then computes Haar scaling and proves the modular-function identity~\eqref{eq:modular}.

Appendix~\ref{sec:counting} gives a second proof of non-unimodularity. It starts with the boundary quasi-action, forms a measure-preserving action using its Radon--Nikodym cocycle, and applies a counting identity for unimodular groups. It uses neither the geometric realization of the boundary closures nor the Haar computation on those realizations. Both proofs use the factorization and boundary estimates from Sections~\ref{sec:models} and~\ref{sec:product}; the appendix isolates a different measure-theoretic route to the obstruction.

An action is \emph{geometric} if it is isometric, continuous, proper and cocompact. We use $(k,c)$-quasi-isometry to mean a map with $c$-dense image satisfying the distance bounds $k^{-1}d-c\le d'\le kd+c$. A uniform quasi-action consists of such maps $q_g$ with $q_e$ uniformly close to the identity and $q_gq_{g'}$ uniformly close to $q_{gg'}$. All bounded errors in statements about an action are uniform in both the point and the group element. We put $a_+=\max(a,0)$ and $\log0=-\infty$.

\section{Model spaces and boundary scales}\label{sec:models}
\subsection{The three models}
The examples differ in the way two upward rays approach. In a smooth factor their distance tends to zero, in a tree they merge, and in a millefeuille space both phenomena occur.
\begin{defn}\label{def:model}
A model space is a pointed space $(X,\omega,h)$ of one of the following types.
\begin{itemize}
\item[(H)] Let $N$ be a nontrivial simply connected nilpotent Lie group and $D$ a derivation of its Lie algebra with eigenvalues of positive real part. Write $e^{tD}$ also for the induced automorphism of $N$. On $N\rtimes_D\R$ take a negatively curved left-invariant metric
\[
           g=(e^{-tD})^*g_N+dt^2,
\]
where $g_N$ is left invariant. The height is $t$, the upper endpoint of the lines $(x,t)$ is $\omega$, and $m_X=\tr D$.
\item[(T)] Let $X$ be a tree of valence $p+1\ge3$, with edges of length $\ell>0$. Choose an end $\omega$. Take $h$ to be minus its Busemann function, increasing towards $\omega$ and taking values in $\ell\Z$ at the vertices. Put $m_X=(\log p)/\ell$.
\item[(M)] Let $M$ have type (H) and $T$ have type (T). Take the equal-height fibre product
\[
             X=\{(u,v)\in M\times T:h_M(u)=h_T(v)\}.
\]
Each vertical line of $T$ specifies a leaf isometric to $M$. Use the length metric obtained by gluing these leaves along their common horoballs \cite[Section~7.A]{CCMT15}. The common height is $h$, the common upper endpoint is $\omega$, and $m_X=m_M+m_T$.
\end{itemize}
\end{defn}
\begin{rem}\label{rem:models}
For any expanding $D$, Heintze's construction gives an inner product and a multiple $\lambda D$ with negative sectional curvature for all sufficiently large $\lambda$ \cite{He74}. The choice of inner product matters. Multiplying all distances and the height of an existing model by $a>0$ divides its exponent by $a$. For (H), use the coordinate $t'=at$, derivation $D/a$ and horizontal metric $a^2g_N$. The tree and mixed assertions follow by rescaling edges and leaves.
\end{rem}
In type (H) the unit gradient lines of $h$ are the vertical geodesics. Left invariance and contraction show that $-h$ is a Busemann function. Homogeneity bounds the sectional curvature above by a negative constant. In type (M), the gluing gives a proper complete $\operatorname{CAT}(-\kappa)$ space for some $\kappa>0$, with $-h$ a Busemann function \cite[Section~7.A]{CCMT15}, \cite[Section~19.2.6]{Cor18}. The projections to $M$ and $T$ are $1$-Lipschitz. Indeed they are so on each leaf, and the gluing distance is computed by chains of segments in leaves.

We will use geodesic completeness in applying Ferragut's theorem. Here is the argument at a branching point of (M). Let a segment end at $q$ of height $z$. Since its height is concave, a sufficiently short terminal subsegment either stays at or above $z$, or stays below $z$ except at $q$. In the tree take a ball about the projection of $q$ meeting no other vertex. The upper part is a single segment and the lower part is a union of disjoint segments. The continuity of the projection therefore puts the terminal subsegment in a single leaf. Extend there, using completeness of $M$. Leaves are convex, so concatenation gives a local geodesic, hence a geodesic \cite[II.1.4]{BH99}. The same argument applies at the other end. Thus all three spaces have the completeness and hyperbolicity required below.

\subsection{Distances seen from a vertical line}
Write $V_\zeta$ for the line with lower endpoint $\zeta$, parametrized by height. The lower endpoints form $\Omega_X=N$ in (H), $\partial T\setminus\{\omega\}$ in (T), and $\Omega_M\times\Omega_T$ in (M). The lines are $(\zeta,t)$, the lines joining the two ends of the tree, and the equal-height pairs of lines, respectively. Every point is on such a line. In (M) this description exhausts the vertical lines because each projection is $1$-Lipschitz and has unit height speed. Evaluation $(\zeta,t)\mapsto V_\zeta(t)$ is continuous. In the tree, nearby lower ends agree above any fixed height once they are sufficiently close, and continuity in a mixed leaf then follows from continuity in $M$. Choose a Borel lower endpoint $\zeta_X(x)$ through $x$ by ordering each vertex's children and following the first one downwards. Use the nilpotent coordinate in (H) and combine the two choices in (M).

For the tree put $\rho_T(\xi,\eta)=\exp h(\xi\wedge\eta)$ for distinct ends, where $\xi\wedge\eta$ is the merging vertex. The shadow of a vertex $v$ is the closed $\rho_T$-ball of radius $e^{h(v)}$ about any of its lower ends and receives mass $e^{m_T h(v)}$. At this radius the open ball is a child shadow and has $1/p$ of that mass. We use closed balls for these exact shadow identities. The Ahlfors estimates below hold for either convention. It splits into $p$ child shadows of equal mass. Finite unions of compact open shadows form a ring and compactness gives countable additivity there, hence a Radon measure $\nu_T$. In the smooth case we need a gauge with the same dilation law.
\begin{lem}\label{lem:gauge}
There is a continuous proper symmetric gauge $r$ on $N$, vanishing only at $e$, such that $r(e^{tD}x)=e^tr(x)$ and $r(xy)\le K\max(r(x),r(y))$. For $\rho_N(x,y)=r(x^{-1}y)$ and Haar measure $\nu_N$, a ball of radius $s$ has measure $\nu_N(B(e,1))s^{\tr D}$.
\end{lem}
\begin{proof}
Fix Euclidean coordinates on the Lie algebra and put $Q=\int_0^\infty e^{-uD^\top}e^{-uD}\,du$. Then $D^\top Q+QD=I$. For $v\ne0$ the derivative of $\|e^{-uD}v\|_Q^2$ is $-\|e^{-uD}v\|^2$, and the two limits as $u$ increases are infinity and zero. If this norm equals one at $u(v)$, set $r(\exp v)=e^{u(v)}$. Its unit ball is the exponential of a compact ellipsoid. This gives continuity, symmetry and homogeneity. Bound $r(xy)$ on the product of two unit balls and dilate to obtain the inequality. Finally $e^{tD}$ has determinant $e^{t\tr D}$.
\end{proof}
Use this gauge and Haar measure in (H). In (M) use $\rho_X=\max(\rho_M,\rho_T)$ and $\nu_X=\nu_M\otimes\nu_T$. Set all gauges to zero on the diagonal.

To see which part of a journey the boundary measures, start with two points of a tree. If the higher point already lies on the upward ray from the lower one, the whole journey costs the difference of their heights. Otherwise one must first go above both points to their merging vertex. Each unit of this extra ascent costs two units of travel, one on each side of the vertex. In a smooth factor the two lines only come close at the corresponding height. We need a bound on the cost of crossing between them that works at every height and on every pair of lines. In the mixed factor the crossing has to be possible in both coordinates at once.

\begin{lem}\label{lem:dist}
The symmetric boundary gauge is continuous, induces the punctured boundary topology, and satisfies $\rho(\xi,\eta)\le K\max(\rho(\xi,\zeta),\rho(\zeta,\eta))$. Its closed balls are compact. There are constants $A\ge1$ and $C_0\ge0$ depending only on the model such that
\[
 A^{-1}r^{m_X}\le\nu_X(B(\zeta,r))\le A r^{m_X}\qquad(\zeta\in\Omega_X,\ r>0).
\]
For $u=\log\rho_X(\xi,\eta)$ and $T=\max(t,t')$ one has
\begin{equation}\label{eq:factor}
 \left|d\big(V_\xi(t),V_\eta(t')\big)-\bigl(|t-t'|+2(u-T)_+\bigr)\right|\le C_0 .
\end{equation}
If $u>T$, every joining geodesic reaches height at least $u-C_0$.
\end{lem}
\begin{proof}
\emph{The tree.} For distinct lower ends, $u$ is the height where the two lines merge. If $u\le T$, the two points are on one vertical line and their distance is $|t-t'|$. If $u>T$, their geodesic goes through the merging vertex, so its length is
\[
                  (u-t)+(u-t')=|t-t'|+2(u-T).
\]
Coincident lower ends give the first case as well. Thus one can take $C_0=0$. For the measure assertion, $B(\zeta,r)$ is the shadow of the vertex on $V_\zeta$ at height $v=\ell\lfloor\log(r)/\ell\rfloor$. Since $v\le\log r<v+\ell$ and $e^{m_T\ell}=p$, its mass lies between $r^{m_T}/p$ and $r^{m_T}$.

\emph{The Heintze space.} The gauge and measure assertions are those of Lemma~\ref{lem:gauge}. Write $P=(x,t)$ and $Q=(y,t')$. Left translation by $P^{-1}$ sends the two points to $(e,0)$ and $(e^{-tD}(x^{-1}y),t'-t)$. The gauge of the new horizontal coordinate is $e^{-t}\rho_X(x,y)$. This gives two uniform comparisons. For each $R\ge0$ choose $M_R\ge1$:
\begin{align}
 d(P,Q)\le R
   &\ \Longrightarrow\ |t-t'|\le R,\quad \rho_X(x,y)\le M_Re^T ,
       \label{eq:local-gauge}\\
 |t-t'|\le R,\quad \rho_X(x,y)\le R e^T
   &\ \Longrightarrow\ d(P,Q)\le M_R .                 \label{eq:local-cross}
\end{align}
For the first implication maximize the continuous gauge on the projection of the compact radius-$R$ ball about $(e,0)$, and use the $1$-Lipschitz height. For the second, the translated point has height in $[-R,R]$ and horizontal gauge at most $R e^{T-t}\le R e^R$. These bounds define a compact set. Increase $M_R$ to cover both comparisons. In particular at any height $s\ge u$, the cost of crossing from $(x,s)$ to $(y,s)$ is at most $M_1$.

Set $s=\max(u,T)$. Ascending to $s$, crossing there, and descending gives
\[
 d(P,Q)\le (s-t)+M_1+(s-t')=|t-t'|+2(u-T)_++M_1.
\]
When $u\le T$ the height difference gives the matching lower bound. Suppose $u>T$. We must show that a different path cannot save a large part of the ascent. Join $(x,u)$ to $(y,u)$ by a segment of length at most $M_1$. Together with the two vertical segments and any geodesic $[P,Q]$ it forms a geodesic quadrilateral. If triangles are $\delta$-thin, each point of $[P,Q]$ is within $R_1=2\delta+M_1$ of one of the vertical segments. The closed $R_1$-neighbourhoods of these segments cover $[P,Q]$, and each contains its corresponding endpoint. Connectedness gives a point $w$ in both neighbourhoods. Choose $(x,s_1)$ and $(y,s_2)$ within $R_1$ of $w$. Their mutual distance is at most $2R_1$, and \eqref{eq:local-gauge} gives
\[
          e^u\le M_{2R_1}e^{\max(s_1,s_2)},\qquad
          h(w)\ge\max(s_1,s_2)-R_1\ge u-E,
          \quad E=\log M_{2R_1}+R_1 .
\]
Since $w$ is on the geodesic, $d(P,Q)\ge (h(w)-t)+(h(w)-t')\ge 2u-t-t'-2E$. This proves both assertions with $C_0=\max(M_1,2E)$. Its independence of $x,y,t,t'$ is what allows us to apply the estimate to a whole family of quasi-isometries below.

\emph{The millefeuille space.} Its boundary balls are products of balls of the same radius. Multiplication of their masses gives the exponent $m_M+m_T$. The maximum gauge has the stated quasi-ultrametric inequality, compact balls and product topology. In each of the three models the lower endpoint topology agrees with the topology of these balls: convergence of endpoints is detected by fellow travel of the corresponding vertical rays down to arbitrarily low fixed heights, as is also seen from \eqref{eq:factor}.

For the distance calculation write $u_M,u_T$ for the two logarithmic boundary separations and $u=\max(u_M,u_T)$. Let $C_M$ be the constant just proved for the smooth factor. The two $1$-Lipschitz projections give
\[
 d(P,Q)\ge\max\{d_M(P_M,Q_M),d_T(P_T,Q_T)\}
       \ge |t-t'|+2(u-T)_+-C_M .
\]
At height $s=\max(u,T)$ the tree coordinates agree. The two points there lie in a common leaf, where their smooth distance is at most $C_M$ since $s\ge u_M$. The path up to $s$, across in the leaf, and down therefore gives the reverse inequality with error $C_M$.

Finally suppose $u>T$. If $u=u_T$, the tree projection of any joining segment contains the merging vertex of height $u$. If $u=u_M$, project the joining segment to the smooth factor and parametrize the projection by arclength. Its length $L$ satisfies
\[
       d_M(P_M,Q_M)\le L\le d(P,Q)\le d_M(P_M,Q_M)+2C_M.
\]
Every subarc has excess length at most $2C_M$ as well, since shortening it would shorten the whole path by the same amount. The projected path is thus a $(1,2C_M)$-quasi-geodesic. By the Morse lemma \cite[III.H.1.7]{BH99} its Hausdorff distance from the smooth joining geodesic is bounded by a constant $R_M$ depending only on $M$. The latter geodesic reaches height $u-C_M$ or higher. The original segment therefore reaches height at least $u-C_M-R_M$. This proves the required uniform height bound also in the mixed case.
\end{proof}

\subsection{Recovering a height scale}
The next question is how a quasi-isometry can change the height coordinate. On a single vertical line even $t\mapsto 2t$ is a quasi-isometry. A model space, however, has many vertical lines, and the height at which two of them come close records the separation of their lower endpoints. We first recover an endpoint from the image of one line. We then use pairs of lines to detect orientation, and families of lines to determine the slope.

We shall use a simple order bound on the real line. A $(k,c)$-quasi-isometry $\phi$ of $\R$ permutes its two ends. If it preserves them, there is $c_1=c_1(k,c)$ with $\phi(\tau)\le\phi(t)+c_1$ whenever $\tau\le t$. In fact, if $\phi(\tau)>\phi(t)$, take unit steps from $\tau$ towards $-\infty$. The values eventually fall below $\phi(t)$, and two successive values differ by at most $k+c$. Hence some $\sigma\le\tau$ satisfies $|\phi(\sigma)-\phi(t)|\le k+c$. The lower quasi-isometry inequality bounds $t-\sigma$ by $k(k+2c)$, and the upper inequality now gives $c_1=k^2(k+2c)+c$.

\begin{lem}\label{lem:line}
Let $P:\R\to X$ be a $(k,c)$-quasi-geodesic. Suppose $|h(P(t))-\phi(t)|\le c$ for an orientation-preserving $(k,c)$-quasi-isometry $\phi$ of the line. There is a unique $\xi\in\Omega_X$ for which $d(P(t),V_\xi(\phi(t)))$ is bounded in $t$. The bound may be chosen in terms of $k,c$ and $X$ alone.
\end{lem}
\begin{proof}
Choose $\xi_t=\zeta_X(P(t))$. These choices can vary from point to point. We shall see that far down the path all of them lie in a single small boundary ball. Let $R$ be a Morse bound for subsegments of $P$, let $C_0$ be as in Lemma~\ref{lem:dist}, and put $B=C_0+R+c+c_1$. For $\tau\le t$ we claim
\[
                     \log\rho(\xi_\tau,\xi_t)\le\phi(t)+B .
\]
If the left side $u$ were larger, it would exceed both endpoint heights by coarse monotonicity. A joining geodesic then reaches at least $u-C_0$, and some $P(\sigma)$, $\tau\le\sigma\le t$, is within $R$ of that high point. Hence $\phi(\sigma)\ge u-C_0-R-c>\phi(t)+c_1$, a contradiction.

It follows from the quasi-ultrametric inequality that all $\xi_\tau,\xi_{\tau'}$ with $\tau,\tau'\le t$ are at distance at most $K e^{\phi(t)+B}$. This tends to zero towards the negative end. Compactness of closed balls yields a unique limit $\xi$, with $\log\rho(\xi_t,\xi)\le\phi(t)+B+\log K$. The two points $P(t)$ and $V_\xi(\phi(t))$ differ in height by at most $c$. Formula~\eqref{eq:factor} therefore bounds their distance by $c+2(B+\log K)+C_0$. If a second lower end had bounded tracking error, its line and $V_\xi$ would be uniformly close at the heights $\phi(t)$. For distinct lower ends \eqref{eq:factor} makes that distance tend to infinity as $t\to-\infty$.
\end{proof}

For two distinct lines in $X$, the number $\log\rho(\xi,\eta)$ is the threshold above which they stay close. A map that takes levels to levels must take such thresholds to thresholds, with a uniform error. This identifies its boundary map. The additional information that fixes its height slope comes from counting how many smaller boundary balls fit in a larger one.
\begin{lem}\label{lem:height}
Suppose $\Psi:X\to X'$ is a quasi-isometry and $h'\Psi=\phi h+O(1)$ for a quasi-isometry $\phi$ of the line. Then $\phi$ preserves orientation. There is a boundary homeomorphism $u$ such that
\begin{align*}
 d'(\Psi V_\xi(t),V'_{u\xi}(\phi(t)))&=O(1),\\
 \log\rho'(u\xi,u\eta)&=\phi(\log\rho(\xi,\eta))+O(1)\quad(\xi\ne\eta),\\
 \phi(s+a)-\phi(s)&=(m_X/m_{X'})a+O(1).
\end{align*}
All errors are uniform in the displayed variables and depend only on the spaces and input bounds. Bounded changes of $\Psi$ preserve $u$.
\end{lem}
\begin{proof}
Enlarge $k,c$ to cover both quasi-isometries and the height error. Apply Lemma~\ref{lem:line} to $t\mapsto\Psi V_\xi(t)$, first replacing $t$ by $-t$ if $\phi$ reverses orientation. With $u\xi$ the resulting lower end, the line estimate has one bound $M$ for all $\xi,t$.

\emph{Detecting the threshold.} Put $s=\log\rho(\xi,\eta)$ and $s'=\log\rho'(u\xi,u\eta)$, provisionally allowing $s'=-\infty$. Let $C_0,C_0'$ be the two distance errors. At heights $t\ge s$ the source lines are at distance at most $C_0$. At $t\le s$ their distance is at least $2(s-t)-C_0$. Applying $\Psi$, moving both image points at most $M$ to their target lines, and using~\eqref{eq:factor} in $X'$ gives
\begin{align}
  s'&\le\phi(t)+\tfrac12(kC_0+c+2M+C_0') &&(t\ge s), \label{eq:upper-threshold}\\
  2(s'-\phi(t))_+&\ge k^{-1}\bigl(2(s-t)-C_0\bigr)-c-2M-C_0' &&(t\le s).
                                                    \label{eq:lower-threshold}
\end{align}
If $\phi$ reversed orientation, sending $t$ to $+\infty$ in~\eqref{eq:upper-threshold} would force $s'=-\infty$, whereas the right side of~\eqref{eq:lower-threshold} becomes positive for sufficiently negative $t$. Thus $\phi$ preserves orientation. To bound $s'$ above, take $t=s$ in~\eqref{eq:upper-threshold}. To bound it below, take $t=s-T_0$ in~\eqref{eq:lower-threshold}, where $T_0>(C_0+k(c+2M+C_0'))/2$ is fixed. Then $s'>\phi(s-T_0)\ge\phi(s)-kT_0-c$. We have proved $|s'-\phi(s)|\le E$ with $E$ independent of the two ends.

A coarse inverse of $\Psi$ respects heights up to bounded error with a coarse inverse of $\phi$ as its height map. Indeed applying $\Psi$ after the coarse inverse moves a point of $X'$ a bounded distance and hence changes its height by a bounded amount. Applying the preceding construction in the reverse direction now gives a boundary map $\bar u$. The two compositions track the original lines at heights within bounded distance of the original parameters. Uniqueness in Lemma~\ref{lem:line} gives $\bar uu=\mathrm{id}$ and $u\bar u=\mathrm{id}$. Since $\phi(s)\to-\infty$ as $s\to-\infty$, the threshold estimate gives continuity of $u$, and the reverse estimate gives continuity of $\bar u$. The same tracking argument shows that maps at bounded distance have the same boundary map.

\emph{Measuring the image of a ball.} Write $B_s(\xi)=B(\xi,e^s)$ and use the analogous notation $B'_v$ in $X'$. Let $D>E+c_1$. The threshold estimate and coarse monotonicity show that
\begin{equation}\label{eq:ball-sandwich}
   B'_{\phi(s)-D}(u\xi)\ \subset\ uB_s(\xi)\
                        \subset\ B'_{\phi(s)+D}(u\xi).
\end{equation}
For the right inclusion, an end in $B_s(\xi)$ has logarithmic separation $r\le s$, so its image has logarithmic separation at most $\phi(r)+E\le\phi(s)+c_1+E$. For the left inclusion, write any other target end as $u\eta$ by surjectivity. If $r=\log\rho(\xi,\eta)>s$, its target separation in logarithmic units is at least $\phi(r)-E\ge\phi(s)-c_1-E>\phi(s)-D$. Thus it cannot lie in the smaller ball. The centre itself belongs to all three sets. With $m'=m_{X'}$, target Ahlfors regularity gives
\begin{equation}\label{eq:ball-image-mass}
            L^{-1}e^{m'\phi(s)}
            \le\nu_{X'}(uB_s(\xi))\le L e^{m'\phi(s)}
\end{equation}
for one $L$ and all $\xi,s$.

\emph{Counting the same centres twice.} Fix $a\ge0$. Within $B_{s+a}(\xi)$ choose a maximal family $I$ of centres whose mutual distances exceed $K e^s$. The balls $B_s(i)$, $i\in I$, are disjoint: an intersection would put their centres at distance at most $K e^s$. Their union lies in $B_{s+a+\log K}(\xi)$. By source Ahlfors regularity these facts bound the size of every finite such family, so one can choose a finite maximal family. Maximality gives the reverse cover
\[
       B_{s+a}(\xi)\ \subset\ \bigcup_{i\in I} B_{s+\log K}(i).
\]
Write $m=m_X$. Apply the source ball bounds to the disjoint union and then to the cover:
\[
 |I|A^{-1}e^{ms}\le A e^{m(s+a+\log K)},\qquad
 A^{-1}e^{m(s+a)}\le |I|A e^{m(s+\log K)} .
\]
Consequently $|I|\asymp e^{ma}$ with constants independent of the centre, $s$ and $a$.

Under the homeomorphism $u$ the small sets remain disjoint, and the enlarged sets still cover the image of the large ball. We can therefore take target measures in exactly the same two inclusions. Using~\eqref{eq:ball-image-mass} gives
\begin{align*}
 |I|L^{-1}e^{m'\phi(s)}
       &\le L e^{m'\phi(s+a+\log K)},\\
 L^{-1}e^{m'\phi(s+a)}
       &\le |I|L e^{m'\phi(s+\log K)} .
\end{align*}
The fixed enlargement by $K$ costs at most $k\log K+c$ in the value of $\phi$, at every level. These inequalities thus give $|I|\asymp e^{m'(\phi(s+a)-\phi(s))}$. Comparison with the source count and passage to logarithms yield
\[
        \big|m'(\phi(s+a)-\phi(s))-ma\big|\le C
\]
for a single $C$. For $a<0$ apply the result starting at $s+a$ with increment $-a$. Dividing by $m'>0$ proves the last assertion.
\end{proof}
For a self quasi-isometry the two exponents coincide. Taking $s=0$ in the last estimate shows that the whole height change is within a fixed distance of the translation $t\mapsto t+\phi(0)$. The constant $\phi(0)$ can vary with the map. The error bound is common to any family with fixed input bounds, the property needed for the group action in Section~\ref{sec:product}.

The affine case controls the full distance, and there is also a converse construction from a single boundary map.
\begin{lem}\label{lem:lift}
If $h'\Psi=\lambda h+c'+O(1)$ with $\lambda>0$ and $\Psi$ a quasi-isometry, then $m_{X'}=m_X/\lambda$ and $d'(\Psi x,\Psi y)=\lambda d(x,y)+O(1)$. If $f$ is a boundary $(L,s)$-quasi-similarity, meaning its scale is $e^s$ up to the factor $L$, put $F_f(x)=V_{f\zeta_X(x)}(h(x)+s)$. This is a rough isometry and
$d(F_f V_\xi(t),V_{f\xi}(t+s))=O(1)$ with bounds depending only on $L,X$.
\end{lem}
\begin{proof}
Insert the line and scale estimates of Lemma~\ref{lem:height} in~\eqref{eq:factor}. Its last estimate identifies $\lambda$ by letting $a\to\infty$. In the converse direction both the merging height and the evaluation height gain $s$, up to $\log L$ for the former. This proves the distance bound. If two lines pass through the same point at height $t$, their log separation is at most $t+O(1)$; applying $f$ shows that either choice of line gives the asserted lift estimate. Surjectivity of $f$ now gives coarse density.
\end{proof}

\section{From a product quasi-isometry to two boundary actions}\label{sec:product}
\subsection{The metric and the examples}
Equip $Z=X\hor Y$ with the length metric associated to a fixed admissible norm on $\R^2$ as in Ferragut's work on geodesics and the visual boundary \cite{Fe26a}. Thus a simultaneous unit vertical motion has unit speed, and the norm is bounded below by $(|a|+|b|)/2$ and above by $C(|a|+|b|)$. The induced topology is the product topology on the closed subset $Z$. Indeed to go from $(x,y)$ to $(x',y')$, follow a geodesic from $x$ to $x'$ while compensating vertically in $Y$. From the resulting point $(x',y_1)$ follow a geodesic from $y_1$ to $y'$ while compensating along a vertical line through $x'$. The second compensation returns to $x'$ since $h_Y(y_1)=h_Y(y')$. The lengths are at most $2C d_X(x,x')$ and $2C(d_X(x,x')+d_Y(y,y'))$. The opposite bound follows from the norm lower bound. Thus $Z$ is complete, locally compact and a length space, hence proper and geodesic \cite[I.3.7]{BH99}.

The continuous surjection $\Theta(\xi,\eta,t)=(V_\xi(t),V_\eta(-t))$ records one lower end in each factor. The same vertical movement counts in both factor distances.
\begin{lem}\label{lem:productdist}
For $p=\Theta(\xi,\eta,t)$ and $q=\Theta(\xi',\eta',t')$,
\begin{align*}
 d_Z(p,q)&=d_X(p_X,q_X)+d_Y(p_Y,q_Y)-|t-t'|+O(1)\\
 &=|t-t'|+2(\log\rho_X(\xi,\xi')-\max(t,t'))_+\\
 &\hspace{34pt}+2(\log\rho_Y(\eta,\eta')+\min(t,t'))_++O(1).
\end{align*}
Each coordinate distance is at most $d_Z(p,q)+O(1)$.
\end{lem}
\begin{proof}
The first equality is \cite[Theorem~1.11]{Fe26b}. Apply~\eqref{eq:factor} twice for the second. Each factor distance dominates $|t-t'|$.
\end{proof}
\begin{ex}\label{ex:models}
(a) In $S=(N_+\times N_-)\rtimes\R$ let $t$ act by $(e^{tD_+},e^{-tD_-})$. Choose $\lambda>0$ so both $\lambda D_\pm$ yield Heintze metrics. Then $(x,y,t)\mapsto((x,t/\lambda),(y,-t/\lambda))$ identifies $S$ with a simply transitive group of product isometries. Its exponents are $\lambda\tr D_\pm$. Left Haar measure in group coordinates is $e^{-t(\tr D_+-\tr D_-)}\,dx\,dy\,dt$, giving $\Delta_S(x,y,t)=e^{-t(\tr D_+-\tr D_-)}$.

(b) Unit edges in trees of valences $p+1,q+1$ give exponents $\log p,\log q$. The product has the Diestel--Leader graph $\DL(p,q)$ as its vertex model \cite[Section~2]{Wo13}.

(c) Take the plane with metric $e^{-2t}dx^2+dt^2$ and a $(p+1)$-regular tree of edge length $\log q$, $q>1$. The exponents are $1,(\log p)/(\log q)$. The geometric action of the same group on this product and on $\HT(q,p)$ identifies their quasi-isometry classes \cite[Theorem~2.15]{BSSW15}. In upper half-plane coordinates $\log y=-h_T$.

(d) A mixed factor can occur in the geometry of a finitely generated group. The following example of Caplinger and Levitin \cite[Section~1.2]{CL25} shows explicitly how branching compensates for unequal smooth expansion rates. Let
\[
 A=\begin{pmatrix}3&1\\1&1\end{pmatrix},\qquad
 \Gamma=\langle \Z^2,s\mid s v s^{-1}=Av\ \ (v\in\Z^2)\rangle .
\]
Here $\Z^2$ denotes the abelian base group. Since $\det A=2$, its image $A\Z^2$ has index two. In the Bass--Serre tree of this ascending HNN extension each vertex has one edge towards the fixed end and two edges in the other direction. This is where the $3$-regular tree enters the geometry.

To identify the smooth directions, diagonalize the real action. The characteristic polynomial is $u^2-4u+2$, with eigenvalues $\alpha=2+\sqrt2>1$ and $\beta=2-\sqrt2\in(0,1)$. Normalize the tree edges to length one and take one application of the stable letter as one unit of absolute height change. On the expanding eigendirection the smooth factor $X$ is the Heintze plane with dilation $x\mapsto\alpha^t x$. On the contracting eigendirection reverse the height, so that the smooth dilation is $y\mapsto\beta^{-t}y$. The second factor $Y$ joins this plane to the tree at equal heights. The geometric action constructed in \cite[Section~1.2]{CL25} is on $X\hor Y$, with opposite heights in these two factors. Their exponents in this normalization are
\[
       m_X=\log\alpha,\qquad
       m_Y=\log(\beta^{-1})+\log 2=\log\alpha ,
\]
where the last equality is $\alpha\beta=\det A=2$. The two smooth exponents alone differ by $\log 2$. The branching term contributes exactly this difference. Thus this discrete geometric action occurs at equality, outside the unequal-exponent hypothesis of Theorem~\ref{thm:A}. A different construction with two identical mixed factors appears in \cite[Section~5]{Cap26}.
\end{ex}
\subsection{Why Ferragut's theorem applies}
The input theorem requires measures on height slices. Their mass on balls in a slice must be uniformly controlled, and the mass of all descendants of a set must grow at the prescribed exponential rate. These requirements explain the choice of exponent in the mixed model.
\begin{lem}\label{lem:admiss}
After multiplying all distances and heights by a common factor $a>0$ chosen to make curvature at most $-1$ in the $\operatorname{CAT}$ sense, the model spaces are horo-pointed admissible. Their exponents in the rescaled height units are $m_X/a$ and $m_Y/a$. Any pair is an admissible couple in the sense of \cite[Definition~3.1 and Remark~3.7]{Fe26b}.
\end{lem}
\begin{proof}
We first compute in the original height units. Replacing $z$ by $z'=az$ transforms $e^{m(z_0-z)}$ into $e^{(m/a)(z'_0-z')}$, giving the rescaled exponents in the statement. Their order is preserved, and $(m_X/a-m_Y/a)(a\chi)=(m_X-m_Y)\chi$. In a Heintze space admissibility is \cite[Section~5.1]{Fe26b}. Its slice measure at $z$ is $e^{-z\tr D}\,dn$ for a fixed Haar measure on $N$. In a tree, edge length measure disintegrates as counting measure on $T_z=h^{-1}(z)$ against $dz$. For a ball of radius $r$ in this slice the count is $p^j$, where $j$ counts vertex heights in $(z,z+r/2]$. Indeed two upward rays must merge by height $z+r/2$. Similarly each point at $z_0\ge z$ has $p^j$ descendants at $z$, now counting vertex heights in $(z,z_0]$. This count lies between $p^{-1}e^{m_T(z_0-z)}$ and $pe^{m_T(z_0-z)}$. Distinct points have disjoint sets of descendants. This gives (E1)--(E3) of \cite{Fe26b}.

For the mixed space use the coordinates $(n,v)\in N\times T$ and the Radon measure $e^{-h_T(v)\tr D}dn\,d\mathrm{length}_T(v)$. Its slice at $z$ is $e^{-z\tr D}dn\otimes\#_{T_z}$. For a Borel set $U$ in the slice at $z_0$, let $U_v=\{n:(n,v)\in U\}$ for $v\in T_{z_0}$. If $j$ counts vertex heights in $(z,z_0]$, descent preserves the smooth coordinate and gives $p^j$ disjoint descendants of each $v$. Therefore, with $\operatorname{Desc}_z U$ denoting all descendants at height $z$,
\[
 \mu_z(\operatorname{Desc}_z U)
 =e^{-z\tr D}p^j\sum_{v\in T_{z_0}}\nu_N(U_v)
 =p^j e^{(z_0-z)\tr D}\mu_{z_0}(U).
\]
The formula holds also for infinite mass, by countable additivity. Since $p^j\asymp e^{m_T(z_0-z)}$, it gives exponent $\tr D+m_T$ for every Borel set. For the slice ball of fixed radius $R>C_0$, with $C_0$ the error in~\eqref{eq:factor}, take boundary product balls of radii $e^{z+(R-C_0)/2}$ and $e^{z+(R+C_0)/2}$. Their evaluations at $z$ sandwich the slice ball. In each case the smooth slice mass is a constant multiple of $e^{(R\pm C_0)\tr D/2}$, and the tree count is bounded above and below uniformly by the preceding calculation. We may choose $R$ at least the threshold $288\delta$ of \cite{Fe26b}.

For completeness, the product of the two slice measures is a Borel kernel in $z$, including when a level contains a branching vertex. Enumerate the countably many open edges $e$ of each tree. For $z$ in the open interval $I_e=h_T(e)$ let $v_e(z)$ be the unique point of $e$ at height $z$. Integration of a nonnegative Borel function over the edge interiors in $T_z$ is the sum over $e$ of its values at $v_e(z)$ multiplied by $\one{I_e}(z)$, each summand extended by zero. Add separately the countable sum over vertices $v$ with indicator $\one{\{h_T(v)\}}(z)$. For a product use these sums in each tree coordinate, at $z$ in $X$ and at $-z$ in $Y$. In smooth coordinates integrate against fixed Haar measures with their continuous height densities. Each term is Borel by parameter integration and the full sum is Borel by monotone convergence. This proves the measurability required for the measure of \cite[Definition~3.6]{Fe26b}, including the exceptional vertex levels.
\end{proof}
An unequal pair of exponents distinguishes the two directions under every quasi-isometry. We order the factors by $m_X>m_Y$ in the proof, interchanging them when necessary.
\begin{cor}[Ferragut's factorization in coordinate form]\label{thm:ferr}
Assume, after relabelling the factors, that $m_X>m_Y$. For every $(k,c)$-quasi-isometry $\Phi$ of $Z$ there are uniform quasi-isometries $\Phi_X,\Phi_Y$ of the factors such that $d_X(\Phi(p)_X,\Phi_X(p_X))$ and $d_Y(\Phi(p)_Y,\Phi_Y(p_Y))$ are bounded uniformly in $p$. All bounds depend only on $k,c$ and the spaces.
\end{cor}
\begin{proof}[Source and uniformity]
This is \cite[Theorem~4.1]{Fe26b} applied to our model spaces via Lemma~\ref{lem:admiss}; we record its coordinate form and uniformity for use below. The quantitative bounds are those in \cite[Section~4.6, (108), (110)]{Fe26b}. The strict inequality between the ordered exponents is a hypothesis of that theorem. Common rescaling preserves the order of exponents. Take the factor maps by holding the opposite lower endpoint fixed. These sections of the product are rough isometric embeddings by Lemma~\ref{lem:productdist}. Composing the sections, $\Phi$, and the projections bounds the coarse Lipschitz constants. The same construction for a coarse inverse of $\Phi$ gives inverse factor maps up to uniform error by applying both coordinate bounds successively.
\end{proof}
\subsection{From bounded error to an exact boundary action}
The factorization separates the coordinates. The height equation makes their scales dependent again.
\begin{lem}\label{lem:oneqi}
There are $L,E$ depending only on the quasi-isometry bounds and the spaces such that each $\Phi$ as above has data $(f,g,b)$ with $f,g$ boundary $(L,b)$- and $(L,-b)$-quasi-similarities and
\begin{equation}\label{eq:oneqi}
          d_Z\big(\Phi\Theta(\xi,\eta,t),\Theta(f\xi,g\eta,t+b)\big)\le E .
\end{equation}
Maps at bounded distance induce the same $f,g$. Their choices of $b$ differ by at most their distance plus the two errors in~\eqref{eq:oneqi}.
\end{lem}
\begin{proof}
At a fixed height, splice $(p_X,p'_Y)$ between $p$ and $p'$. The factorization bounds the change of image height in each of the two comparisons by twice its coordinate error. Thus $h\Phi=\phi h+O(1)$, where $\phi$ is obtained along one fixed vertical line. This function is coarsely Lipschitz, and a coarse inverse of $\Phi$ gives its coarse inverse by the same construction. On the two factors use $\phi(t)$ and $-\phi(-t)$ in Lemma~\ref{lem:height}. Since a factor is compared with itself, its slope is one, giving $\phi(t)=t+b+O(1)$. The line and scale estimates give $f,g$ and the product distance estimate gives~\eqref{eq:oneqi}. For uniqueness compare the two right-hand models for all $t$. Height bounds the shift difference. Sending $t$ to $-\infty$ and $+\infty$ successively in Lemma~\ref{lem:productdist} forces equality of the two lower endpoint maps.
\end{proof}
\begin{prop}\label{prop:qa}
For a uniform quasi-action $(q_a)$ use $(f_a,g_a,b_a)$ from Lemma~\ref{lem:oneqi}. The maps $a\mapsto f_a$ and $a\mapsto g_a$ are homomorphisms, and $|b_{aa'}-b_a-b_{a'}|$ is uniformly bounded. The lower endpoints can be recovered by
\begin{align*}
 f_a\xi&=\lim_{n\to\infty}\zeta_X\big((q_a\Theta(\xi,\eta,-n))_X\big),\\
 g_a\eta&=\lim_{n\to\infty}\zeta_Y\big((q_a\Theta(\xi,\eta,n))_Y\big).
\end{align*}
\end{prop}
\begin{proof}
Applying~\eqref{eq:oneqi} twice gives data $(f_af_{a'},g_ag_{a'},b_a+b_{a'})$ for the composition with a uniform error. Compare with $q_{aa'}$ and apply uniqueness, also at the identity. For the limits, if $x$ is within $R$ of $V_\zeta(s)$ then~\eqref{eq:factor} gives $\log\rho_X(\zeta_X(x),\zeta)\le s+O_R(1)$. In the first formula $s=-n+b_a$, in the second $s=-n-b_a$.
\end{proof}
\begin{lem}\label{lem:sc}
Every compactly generated locally compact $H$ has a compact normal subgroup $K$ with second countable quotient; passage to this quotient preserves the quasi-isometry class, unimodularity and amenability. For second countable $H$ the transferred quasi-action and the two boundary actions can be taken jointly Borel. One may use $b_a=h(q_a o)-h(o)$ in~\eqref{eq:oneqi} with enlarged uniform constants.
\end{lem}
\begin{proof}
Use Kakutani--Kodaira \cite[Theorem~8.7]{HR79}. The quotient by a compact kernel is a quasi-isometry for compact word metrics \cite[Chapter~4]{CdlH16}. Its normalized Haar measure on the kernel gives $\Delta_H=\Delta_{H/K}\circ\mathrm{pr}$; amenability also passes in both directions. Choose quasi-inverses constant on countable Borel partitions into sets of bounded diameter. In $H$ one obtains such a partition from translates of an open $V$ with $V^{-1}V$ in the compact generating neighbourhood. In $Z$ use a countable cover by unit balls. Here is the uniform bound behind these replacements. If the original maps $\varphi,\bar\varphi$ have coarse Lipschitz bounds $(k,c)$ and $(\bar k,\bar c)$, and the partition atoms have diameters at most $\delta_H,\delta_Z$, the replacement maps $\varphi_0,\bar\varphi_0$ satisfy
\[
 d_Z(\varphi_0(x),\varphi(x))\le A:=k\delta_H+c,\qquad
 d_H(\bar\varphi_0(z),\bar\varphi(z))\le A':=\bar k\delta_Z+\bar c.
\]
Left invariance of the word metric then gives, uniformly in $a,z$,
\[
 d_Z(\varphi_0 L_a\bar\varphi_0(z),\varphi L_a\bar\varphi(z))
 \le A+kA'+c.
\]
Thus the same boundary maps and bounded-error conclusions apply after replacement. Renaming these Borel maps, $q_a=\varphi L_a\bar\varphi$ is jointly Borel; the limits of Proposition~\ref{prop:qa} prove the assertions for the boundary. Evaluating~\eqref{eq:oneqi} at $o$ gives the bound on the change of $b_a$.
\end{proof}

\section{Locally compact closures of boundary actions}\label{sec:closure}
There are two obstructions to compactness of a family of similarities of the line: the image of the origin can escape, or the scale can tend to zero or infinity. The same two tests work on the punctured boundary of a model space. Write $G_X$ for the stabilizer of infinity in the homeomorphism group of $\widehat\Omega_X=\Omega_X\cup\{\infty\}$, with its compact-open topology. This is a Polish group \cite[Section~9.B]{Ke95}. Fix distinct $\xi_*,\eta_*$ and set
\[
  \beta(f)=\log\frac{\rho_X(f\xi_*,f\eta_*)}{\rho_X(\xi_*,\eta_*)}.
\]
An $(L,s)$-quasi-similarity has $|\beta(f)-s|\le\log L$ and hence has distortion at most $L^2$ at scale $e^{\beta(f)}$.

\begin{prop}\label{prop:closure}
Let $\Gamma\le G_X$ consist of $L$-quasi-similarities and put $Q=\overline{\Gamma}$. Fix $o\in X$ with chosen lower end $\xi_o=\zeta_X(o)$.
\begin{enumerate}
\item $Q$ is locally compact and second countable, its elements are $L$-quasi-similarities, and its action on the boundary and the function $\beta$ are continuous. Bounds on $|\beta(f)|$ and $\rho_X(f\xi_o,\xi_o)$ define compact closed subsets of $Q$.
\item The maps $F_f(x)=V_{f\zeta_X(x)}(h(x)+\beta(f))$ give a uniform quasi-action by rough isometries. Each orbit map $f\mapsto F_f(x)$ is continuous, and $\{f:d(F_f(o),o)\le R\}$ is compact.
\item If the orbit of $o$ under the lifts of $\Gamma$ is cobounded, then $Q$ is compactly generated and $f\mapsto F_f(o)$ is a quasi-isometry. The group is focal hyperbolic, amenable, and has connected, totally disconnected or mixed type according as $X$ has type (H), (T) or (M).
\end{enumerate}
\end{prop}
\begin{proof}
For (1), if $f_n$ converges to a homeomorphism in $G_X$, evaluation on the reference pair bounds the possible logarithmic scales. Pass to a convergent subsequence of scales and take limits in the inequalities, retaining $L$. The action and $\beta$ are continuous by evaluation. We give the compactness argument with uniform constants. Suppose $|\beta(f)|\le B$ and $\rho_X(f\xi_o,\xi_o)\le R$, and set $C=L^2e^B$. The scale estimate preceding the proposition gives
\[
 C^{-1}\rho_X(\xi,\eta)\le\rho_X(f\xi,f\eta)\le C\rho_X(\xi,\eta).
\]
The same bounds hold for $f^{-1}$, whose displacement at $\xi_o$ is at most $CR$. The quasi-ultrametric inequality therefore puts the images under both $f$ and $f^{-1}$ of every fixed compact ball in a common compact ball. It also gives equicontinuity on these balls. Arzel\`a--Ascoli on a countable exhaustion, applied simultaneously to the maps and their inverses, produces locally uniform limits $f$ and $j$. Equicontinuity and the compact image bounds allow passage to the limit in the two composition identities, so $fj=jf=\mathrm{id}$.

We must also control the point at infinity. If $\rho_X(f\xi,\xi_o)\le S$, then
\[
 \rho_X(\xi,\xi_o)
 \le C\rho_X(f\xi,f\xi_o)
 \le CK\max\{S,R\}.
\]
For $f^{-1}$ the same argument gives the bound $CK\max\{S,CR\}$. Thus, for each $S$, points outside a single sufficiently large ball are sent outside $B(\xi_o,S)$ by every map in the family and by every inverse. The family is equicontinuous at infinity in the one-point compactification. Together with local uniform convergence, this gives uniform convergence there for the maps and their inverses. The constraints are closed, so the constrained set is compact. Such a set contains an identity neighbourhood; second countability follows from that of $G_X$.

For (2) apply Lemma~\ref{lem:lift} with logarithmic scale $\beta(f)$. On a group of $L$-quasi-similarities the defect of $\beta$ is at most $4\log L$, by evaluating the composite on the fixed pair. Thus $F_fF_{f'}$ and $F_{ff'}$ both follow the same line with heights at bounded distance, uniformly; $F_e$ is the identity. Evaluation of a line at a height gives continuity of each orbit map. If $d(F_f(o),o)\le R$, the height change is bounded by $R$ and~\eqref{eq:factor} bounds $\rho_X(f\xi_o,\xi_o)$. The set is closed and hence compact by (1).

For (3) let $R_0$ bound the distance of points in $X$ from the $Q$-orbit. Along a geodesic from $o$ to $F_f(o)$ take points at successive distances at most one. Approximate each by an orbit point, using $e,f$ at the two endpoints. If $F_a(o),F_b(o)$ are successive choices, applying $F_{a^{-1}}$ shows $d(F_{a^{-1}b}(o),o)\le 2R_0+1+C$ for a uniform $C$. Thus the compact symmetric enlargement of this bounded-displacement set generates $Q$, with word length at most $d(F_f(o),o)+2$. Enlarge it to contain a compact identity neighbourhood. Conversely bounded generator displacement and the uniform quasi-action inequalities bound displacement by a linear function of word length: successive orbit points along a word differ by a uniform amount since every $F_a$ has additive distortion bounded independently of $a$. The same comparison applied to $a^{-1}b$ proves the orbit map is a quasi-isometry.

This map identifies $\partial Q$ with $\partial X$, equivariantly for the boundary extensions of $F_f$. Their upward vertical rays have common endpoint $\omega$, fixed by every $f$. Distinct lower ends give distinct boundary points by~\eqref{eq:factor}, so the hyperbolic group $Q$ is non-elementary. The fixed-point characterization gives focality and amenability \cite[Theorem~7.3]{CCMT15}. For a tree the compactification is totally disconnected and stabilizers of finite clopen partitions give open subgroups separating any nonidentity homeomorphism. Thus $Q$ is totally disconnected. For a Heintze space the boundary is a positive-dimensional sphere. Of the three focal types only the connected type has such a boundary: the other types have Cantor boundary or a boundary that is not locally connected \cite[Propositions~19.9, 19.10, 19.15]{Cor18}. In type (M) the model is, after edge normalization, a pure millefeuille space. Its boundary identifies the mixed type \cite[Corollary~19.17 and Proposition~19.21]{Cor18}.
\end{proof}

We now relate the orbit hypothesis to the condition on pairs used in Theorem~\ref{thm:G}. For $P_X=\{(\xi,\eta):\xi\ne\eta\}$ set $\varpi(\xi,\eta)=V_\xi(\log\rho_X(\xi,\eta))$.
\begin{lem}\label{lem:pairs}
The map $\varpi$ is continuous, has coarsely dense image, and the inverse image of each bounded set lies in a compact subset of $P_X$. Moreover $d(F_f\varpi(p),\varpi(fp))$ is uniformly bounded. Thus cocompactness of $\Gamma$ on pairs gives coboundedness of its lifted orbit. Conversely coboundedness of the $Q$-orbit implies that $Q$ acts cocompactly on pairs.
\end{lem}
\begin{proof}
In a tree all vertices are merging vertices. In a smooth factor $\rho_N(\xi,\cdot)$ takes every positive value by dilation; in a mixed factor hold the tree endpoint fixed. These facts give coarse density. If $d(\varpi(\xi,\eta),o)\le R$, the height bounds $\rho_X(\xi,\eta)$ above and away from zero. Equation~\eqref{eq:factor} bounds $\rho_X(\xi,\xi_o)$ above, and the quasi-ultrametric inequality bounds $\rho_X(\eta,\xi_o)$. These closed constraints define a compact subset disjoint from the diagonal. The lift estimate in Lemma~\ref{lem:lift} and the bound $|\log\rho_X(f\xi,f\eta)-\log\rho_X(\xi,\eta)-\beta(f)|\le2\log L$ give coarse equivariance. A compact set of pair representatives now gives a bounded set whose lifted translates cover $X$ to bounded error. Conversely move $\varpi(p)$ close to $F_f(o)$ and apply $F_{f^{-1}}$ to see that $\varpi(f^{-1}p)$ stays in a fixed ball. The compact constraint above gives representatives for all pairs.
\end{proof}
This lemma and Proposition~\ref{prop:closure} prove Theorem~\ref{thm:G}.

\subsection{What the two projections retain}
A single boundary action can have a large kernel. Taken together the two actions retain every group element up to a compact kernel because together they control displacement in the product.
\begin{thm}\label{thm:projections}
Under the hypotheses of Theorem~\ref{thm:A}, suppose $H$ is second countable. Let $Q_X,Q_Y$ be the closures of the boundary images of Proposition~\ref{prop:qa}. The homomorphisms $\pi_X(a)=f_a$ and $\pi_Y(a)=g_a$ are continuous, and both lifted orbits are cobounded. Uniformly in $a,p$,
\[
 \beta_X(f_a)=b_a+O(1),\qquad \beta_Y(g_a)=-b_a+O(1),
\]
and the coordinates of $q_a(p)$ are within bounded distance of the respective lifts $F_{f_a}(p_X),F_{g_a}(p_Y)$. The combined map $\pi:H\to Q_X\times Q_Y$ is proper with closed image and compact kernel. In particular $H$ is amenable.
\end{thm}
\begin{proof}
Use the Borel version of Lemma~\ref{lem:sc}. Evaluations on a countable dense set generate the Borel structure of $G_X$, so $\pi_X$ is a Borel homomorphism of Polish groups and therefore continuous \cite[Theorem~9.10]{Ke95}. The same holds in $Y$. Evaluation on one pair gives the scale bounds. Equations~\eqref{eq:oneqi} and~\eqref{eq:factor} give the coordinate bounds. The orbit $q_a(o)=\varphi(a\bar\varphi(o))$ is cobounded in $Z$, whose two coordinate projections are onto. The projected lifted orbits are consequently cobounded.

For compact $C\subset Q_X\times Q_Y$, continuity of the lifted orbit maps bounds both coordinates of $q_a(o)$ when $\pi(a)\in C$. Lemma~\ref{lem:productdist} bounds $d_Z(q_a(o),o)$. Since $\varphi$ is a quasi-isometry, $a\bar\varphi(o)$, hence $a$, lies in a relatively compact set (compact word balls suffice before a fixed right translation). The preimage of $C$ is closed, proving properness. A proper continuous homomorphism here has closed image, compact kernel and induces a topological isomorphism from the quotient onto the image. Both closures are amenable by Proposition~\ref{prop:closure}. Amenability passes to products, closed subgroups and compact extensions \cite[Section~4.1]{Zi84}.
\end{proof}

\section{Realizing the actions with a common height unit}\label{sec:actions}
A focal group has a geometric action on a negatively curved model. To combine two such actions we must identify both the direction and the unit of their Busemann heights. We first record the geometric facts that make this comparison possible.

\subsection{Geometric models for the boundary closures}
\begin{lem}\label{lem:normalform}
Let $M$ be a simply connected homogeneous negatively curved Riemannian manifold of dimension at least two. Suppose a closed cocompact group of isometries fixes exactly one boundary point $\omega$. For every height $h$ equal to minus a Busemann function at $\omega$, $(M,\omega,h)$ is isometric to a Heintze model with its stated height.
\end{lem}
\begin{proof}
By \cite{He74} a solvable group $S$ acts simply transitively on $M$. It is amenable and non-elementary hyperbolic and fixes a unique boundary point \cite[Theorem~7.3]{CCMT15}. If the full isometry group fixes a point, it must be the point fixed by both the given group and $S$. Otherwise $M$ is symmetric and the full group acts transitively on the boundary \cite[Lemma~19.11]{Cor18}; conjugate $S$ so that it fixes $\omega$.

Identify $M$ with $S$ at a point of height zero. Height increments give a character on $S$, and $X_0=\nabla h(e)$ defines its left-invariant unit gradient field \cite{HIH77}. The kernel $N$ acts simply transitively on the zero horosphere. The map $(x,t)\mapsto x\exp(tX_0)$ identifies $N\times\R$ diffeomorphically with $S$, so $N$ is connected and simply connected. Orthogonality gives the model metric for $D=\operatorname{ad}(X_0)|_{\mathfrak n}$. To see expansion, choose $-a^2<0$ bounding curvature above. The variation $\exp(sv)\exp(tX_0)$ for $v\in\mathfrak n$ gives a stable orthogonal Jacobi field of norm $\|e^{-tD}v\|$, bounded by $e^{-at}\|v\|$ for $t\ge0$ \cite{HIH77}. All eigenvalues of $D$ therefore have positive real part. Generalized eigenspaces satisfy $[\mathfrak n_\alpha,\mathfrak n_\beta]\subset\mathfrak n_{\alpha+\beta}$. Real parts add, so sufficiently long iterated brackets vanish. This proves nilpotence and identifies the height with $t$.
\end{proof}
\begin{lem}\label{lem:focalmodel}
A focal group of connected, totally disconnected or mixed type has a continuous proper cocompact isometric action, fixing the upper end, on a model space of type (H), (T) or (M), respectively.
\end{lem}
\begin{proof}
Use the three geometric actions in \cite[Propositions~19.8--19.10]{Cor18}, following \cite[Theorem~7.3]{CCMT15}. In the connected type the image is closed by properness, and the boundary identification gives exactly one fixed point. Lemma~\ref{lem:normalform} identifies the metric and the height. In the totally disconnected type use the regular tree, with unit edges and its invariant end. In the mixed type use the pure millefeuille realization of \cite[Proposition~19.10]{Cor18}; the boundary characterization is \cite[Proposition~19.21]{Cor18}. Its upper end is fixed by every boundary homeomorphism \cite[Corollary~19.17]{Cor18}. The distinguished smooth end has a transitive stabilizer. Indeed, in the nonsymmetric case the full isometry group fixes its distinguished end and is transitive. In the rank-one symmetric case, conjugate a simply transitive Iwasawa subgroup so that it fixes the prescribed end. In either case the stabilizer is focal by its geometric action and boundary fixed point. The isometric equivalence between the prescribed Busemann functions is supplied by \cite[Lemma~19.13]{Cor18}. The smooth leaf has dimension $\dim N+1\ge2$, since the smooth nilpotent factor is nontrivial. Lemma~\ref{lem:normalform} therefore puts the leaf in Heintze form. Use increasing heights in both the leaf and the tree; changing both signs keeps their equal-height equation.
\end{proof}
An end-fixing isometry has a height increment independent of the point. Its absolute value can be detected without coordinates by following a long orbit.
\begin{lem}\label{lem:stable}
On $A_X=\Isom(X)_\omega$ the formula $\tau(a)=h(ax)-h(x)$ defines a continuous homomorphism, and $\lim_{n\to\infty}d(a^nx,x)/n=|\tau(a)|$.
\end{lem}
\begin{proof}
Busemann functions at the same end differ by a constant \cite[II.8]{BH99}. Evaluation at $x$ proves continuity and composition proves additivity. The height bound gives one inequality. If the stable length is positive, $a$ is loxodromic \cite[Section~8.1]{Gr87}. The two boundary fixed points determine a geodesic line, unique in a proper $\operatorname{CAT}(-\kappa)$ space \cite[II.2.13 and II.9]{BH99}. On that line $a$ translates by $\tau(a)$ in the increasing height parameter. Stable length is independent of the starting point, proving the other inequality. If the stable length is zero, the first inequality gives equality.
\end{proof}

\subsection{Comparing coarse height with an exact character}
An orbit quasi-isometry between two models need not visibly respect height. The next argument first recovers the height change of powers, then deduces a single scale for all elements.
\begin{lem}\label{lem:compare}
Suppose $Q$ has a cobounded uniform quasi-action $(F_a)$ on $X$ with $h(F_a x)=h(x)+\beta(a)+O(1)$, and an isometric action fixing the upper end of another model $U$, with character $\tau$. Let $\Psi:X\to U$ be a quasi-isometry with $d(\Psi F_a x,a\Psi x)$ uniformly bounded. There are constants $\lambda>0,c$ such that
\[
 \tau(a)=\lambda\beta(a)+O(1),\qquad h_U(\Psi x)=\lambda h_X(x)+c+O(1).
\]
The map is a rough similarity of ratio $\lambda$ and $m_U=m_X/\lambda$.
\end{lem}
\begin{proof}
The height estimate and composition law give $\beta$ bounded additive defect. Its homogenization $B(a)=\lim_n\beta(a^n)/n$ exists, has bounded defect $D$, and differs from $\beta$ by a bounded function \cite[Section~2.2.2]{Ca09}. Applying the hypotheses to $a^n$ gives
\[
 |\beta(a^n)|\le d(F_{a^n}o,o)+O(1)\le k\,d(a^n\Psi o,\Psi o)+O(1).
\]
Divide by $n$ and apply Lemma~\ref{lem:stable}: $|B(a)|\le k|\tau(a)|$.

If $\tau(a)=\tau(a')$, then $a^{-n}(a')^n$ has zero character and hence zero $B$ for every $n$. The defect estimate gives $n|B(a)-B(a')|\le D$. Thus $B$ factors \emph{exactly} through a homogeneous quasimorphism on $\tau(Q)\le\R$. On this abelian group it is additive: apply the defect inequality to $nu,nv$ and divide by $n$. The bound by $k|u|$ makes the resulting additive function Lipschitz. It is linear on a cyclic subgroup, and on a dense subgroup extends continuously to $\R$. Thus $B=\mu\tau$. Coboundedness and the unboundedness of height imply $\beta$ is unbounded, so $\mu\ne0$.

Put $\lambda=1/\mu$. The first estimate follows. To obtain the second, approximate $x$ within fixed distance by $F_a o$. Then $h_X(x)=h_X(o)+\beta(a)+O(1)$ and $h_U(\Psi x)=h_U(\Psi o)+\tau(a)+O(1)$. Lemma~\ref{lem:height} forces the resulting affine slope to be positive. Lemma~\ref{lem:lift} proves the distance and exponent assertions.
\end{proof}

\subsection{Joining the factor actions}
\begin{thm}\label{thm:joining}
In the second countable setting of Theorem~\ref{thm:projections}, the geometric model and all three rough isometries in Theorem~\ref{thm:A}(2) exist. The height character on its first factor satisfies $\chi(a)=b_a+O(1)$ uniformly on $H$.
\end{thm}
\begin{proof}
\emph{Choose the two units.} For $Q_X$ take a model $U$ from Lemma~\ref{lem:focalmodel}. The orbit maps $a\mapsto F_a o$ and $a\mapsto a o'$ are quasi-isometries to $X$ and $U$, by Proposition~\ref{prop:closure} and the locally compact Milnor--\v{S}varc lemma \cite[Theorem~4.C.5]{CdlH16}. For each $x$ choose $a_x$ with $F_{a_x}o$ within fixed distance of $x$, and put $\Psi_X(x)=a_x o'$. This is a quasi-isometry. The points $F_{a_{F_ax}}o$ and $F_{aa_x}o$ are a bounded distance apart, uniformly in $a,x$. The orbit comparison bounds the word distance between the two group elements, hence bounds their displacements relative to each other at $o'$. This proves coarse equivariance. Construct $\Psi_Y$ in the same way. Lemma~\ref{lem:compare} gives slopes $\lambda_X,\lambda_Y$. Divide the metric and the height of each target factor by its \emph{own} slope. By Remark~\ref{rem:models} the resulting spaces $X',Y'$ have exponents $m_X,m_Y$ and each factor map is a rough isometry. Write $\tau_X,\tau_Y$ for the rescaled exact characters.

\emph{Make the shifts cancel.} On $H$, Theorem~\ref{thm:projections} and Lemma~\ref{lem:compare} give
\[
 \tau_X(f_a)=b_a+O(1),\qquad \tau_Y(g_a)=-b_a+O(1).
\]
Their sum is a bounded homomorphism to $\R$. Evaluation at powers forces it to vanish identically. The product action therefore preserves $Z'=X'\hor Y'$. Put $\chi(a)=\tau_X(f_a)$; it is continuous and $\chi(a)=\lim_n b_{a^n}/n$.

\emph{Compare the two products.} For $(x,y)\in Z$, the sum of heights of $(\Psi_X x,\Psi_Y y)$ is bounded, since both normalized height maps differ from the original heights by bounded amounts. Move the second point on a chosen vertical line by precisely the amount required to make the sum zero. Denote the resulting map by $F$. The change is uniformly bounded. The first identity in Lemma~\ref{lem:productdist}, applied in both products, now gives
\[
                   d_{Z'}(Fp,Fq)=d_Z(p,q)+O(1).
\]
For coarse density choose preimages close in each target factor. Their original heights have bounded sum; a bounded vertical adjustment in the source and the rough-isometry bounds put an image uniformly close to the prescribed target point. The coordinate estimates of Theorem~\ref{thm:projections}, coarse equivariance on the factors and the bounded correction give $d_{Z'}(Fq_a p,aFp)=O(1)$.

\emph{Check the group topology.} The action is continuous since both factor actions are. Coarse equivariance and coarse density make the orbit cobounded, and properness of $Z'$ makes it cocompact. For a compact $C\subset Z'$, the elements moving $C$ to meet itself have both boundary projections in compact sets, by properness of the individual geometric actions. The combined homomorphism of Theorem~\ref{thm:projections} is proper; the transporter of $C$ is closed, hence compact. The action kernel $W$ is compact and normal, and the induced action is faithful and geometric.
\end{proof}

\section{Haar measure and the exact modular character}\label{sec:mod}
The action constructed in Section~\ref{sec:actions} makes the two height shifts cancel. Their effects on Haar measure need not cancel, since each unit of height carries the exponent of its own factor. We first compute this effect for one factor, then justify its transfer to a group which acts only cocompactly on the product.

\begin{lem}\label{lem:restriction}
Let $D$ be a closed subgroup of a locally compact group $A$.
\begin{enumerate}
\item If $D$ is normal, or if $A$ is amenable and $A/D$ is compact, then $\Delta_A|_D=\Delta_D$.
\item If $D$ is cocompact, two continuous real characters on $A$ that agree on $D$ agree on $A$.
\end{enumerate}
\end{lem}
\begin{proof}
The invariant-measure criterion on a homogeneous space says that $A/D$ carries a nonzero $A$-invariant Radon measure precisely when the two modular functions agree on $D$ \cite[Theorem~2.51]{Fo16}. For a normal subgroup, take left Haar measure on the quotient group. For an amenable group acting on the compact space $A/D$, apply its fixed-point property to the compact convex space of probability measures, giving an invariant probability. For the second assertion, the difference descends to a continuous real function on compact $A/D$. A bounded homomorphism to $\R$ vanishes by evaluation on powers.
\end{proof}
\begin{prop}\label{prop:factorhaar}
Let $U$ be any model space and $A_U=\Isom(U)_{\omega_U}$. Write $\tau_U(a)=h_U(au)-h_U(u)$. Then $A_U$ is locally compact and amenable, and
\begin{equation}\label{eq:factorhaar}
                       \log\Delta_{A_U}(a)=-m_U\tau_U(a).
\end{equation}
\end{prop}
\begin{proof}
We exhibit a closed geometric subgroup in each end stabilizer and do the computation there first. The stabilizer itself is closed in the isometry group of a proper metric space and acts properly. A geometric subgroup makes its action cocompact. Its orbit map then identifies its hyperbolic boundary with $\partial U$ equivariantly. This boundary has more than two points and the group fixes its upper end, so the stabilizer is focal and amenable \cite[Theorem~7.3]{CCMT15}. Also, any closed geometric subgroup $D$ is cocompact in $A_U$: if $U=DC$ with $C$ compact, then $A_U=D\{a:ao\in C\}$ for $o\in U$, and the set in braces is compact by properness. One can therefore pass from $D$ to $A_U$ using both assertions of Lemma~\ref{lem:restriction}.

\emph{Heintze space.} Take $D=N\rtimes_B\R$ acting simply transitively by left translations. Its image is closed because its orbit map is a homeomorphism. In group coordinates left Haar measure is $e^{-s\tr B}dn\,ds$. Right multiplication by $(e,t)$ changes its mass by $e^{-t\tr B}$. Right multiplication by a horizontal translation preserves it since $N$ is unimodular. The height character is $t$ and $m_U=\tr B$, giving~\eqref{eq:factorhaar} on $D$ and hence on $A_U$.

\emph{Regular tree.} The full end stabilizer is transitive on vertices. Let $v^+$ be the parent of $v$ and $K_v$ its compact open stabilizer inside $A_U$. The group $K_{v^+}$ permutes the $p$ children transitively and $[K_{v^+}:K_v]=p$. Choose $s$ translating $v$ to $v^+$ towards the distinguished end on an axis. Then $sK_vs^{-1}=K_{v^+}$. Under our convention conjugation by $s$ multiplies left Haar measure by $\Delta(s)^{-1}$. Thus $\Delta(s)=p^{-1}$ while $\tau_U(s)=\ell$. If $\tau_U(a)=0$, the rays from $v$ and $av$ merge, and the points on them at an equal sufficiently high vertex level coincide. Thus $a$ fixes a vertex and belongs to a compact subgroup, where $\Delta=1$. Since $\tau_U(A_U)=\ell\Z$, the formula follows.

\emph{Millefeuille space.} Write $U=M\times_h T$ for the equal-height fibre product. Take $S=N\rtimes_B\R$ and the tree end stabilizer $A_T$. The group
\[
       D=\{(s,a)\in S\times A_T:\tau_M(s)-\tau_T(a)=0\}
\]
is a closed normal subgroup. It acts on $U$ by leaf isometries compatible with the gluing metric. This action is faithful since both coordinate projections of $U$ are onto. It is proper since compact subsets project to compact sets in the factors and both factor actions are proper. It is cocompact: at each vertex height the tree stabilizer is transitive on vertices at that level, $N$ is transitive on the smooth horosphere, and a matched shift moves between consecutive vertex levels. Every point is within $\ell$ of such a level on a vertical line. Properness and faithfulness identify $D$ with a closed subgroup of $A_U$. By normality its modular function is the restriction of the product modular function. The two shifts now have the \emph{same} sign and common value $\tau_U$, giving $-(m_M+m_T)\tau_U$. Since $m_U=m_M+m_T$, the transfer in the opening paragraph proves~\eqref{eq:factorhaar} on all of $A_U$.
\end{proof}

We turn to the product with \emph{opposite} heights. The full isometry group is unnecessary here. Working in the group of pairs that preserve the defining equation gives access to both factor characters and a proper action.
\begin{thm}\label{thm:modular}
Suppose a locally compact group $G$ acts faithfully and geometrically on $U\hor V$ through a continuous homomorphism into $A_U\times A_V$. If the first height shift is $\tau(g)$ and the second is $-\tau(g)$, then
\[
                       \log\Delta_G(g)=-(m_U-m_V)\tau(g).
\]
\end{thm}
\begin{proof}
The group of height-compatible pairs
\[
              P=\ker(\tau_U+\tau_V:A_U\times A_V\to\R)
\]
is closed, normal and amenable. By Lemma~\ref{lem:restriction} its modular function is the restriction of $\Delta_{A_U}\Delta_{A_V}$. Its action on $U\hor V$ is proper by the two coordinate projections. The map $G\to P$ is injective by faithfulness and proper: a compact set of pairs moves a fixed product point through a compact set, whose inverse image under the $G$-orbit map is compact. Thus it identifies $G$ topologically with a closed subgroup.

To see cocompactness in $P$, choose a compact $C$ with $GC=U\hor V$. For each $p\in P$, write $po=gz$ with $z\in C$. Then $g^{-1}p$ is in the compact set $\{a\in P:ao\in C\}$. Thus $G\backslash P$ and, by inversion, $P/G$ are compact. Apply Lemma~\ref{lem:restriction} to $G\le P$. On its image the two shifts are opposite, so Proposition~\ref{prop:factorhaar} yields
\[
   \log\Delta_G(g)= -m_U\tau_U(g_U)-m_V\tau_V(g_V)
                  =-(m_U-m_V)\tau(g).
\]
\end{proof}
\begin{proof}[Proof of Theorem~\ref{thm:A}]
In the second countable case use Theorem~\ref{thm:joining}. The exponent normalizations and the opposite shifts give~\eqref{eq:modular} on $H/W$ by Theorem~\ref{thm:modular}, and hence on $H$ because $W$ is compact. The character is nonzero: height is unbounded on the product, while a zero character would make it bounded on all translates of a compact fundamental set. The two exponents differ, so $\Delta_H$ is nontrivial. Amenability follows from Theorem~\ref{thm:projections}.

For general $H$, take $K$ as in Lemma~\ref{lem:sc} and pass to $\bar H=H/K$. To track an initially prescribed quasi-isometry, choose a section $\sigma:\bar H\to H$ and use $\varphi\sigma$ on the quotient. Each fibre of $H\to\bar H$ has uniformly bounded word diameter because $K$ is compact and the word metric is left invariant. With coarse inverse $\mathrm{pr}\,\bar\varphi$, the transferred maps for $\mathrm{pr}(a)$ are therefore within bounded distance of $q_a$, uniformly in $a$. Borel replacements on the second countable quotient have the same property. This retains the comparison of $b_a$ and the coarse equivariance of $F$. Lift the compact kernel on $\bar H$ to its compact normal preimage in $H$. Both the modular function and the geometric action pull back through the quotient. Finally, if the factors were interchanged to apply Ferragut's theorem, restore their order; this reverses both signs in the product $(m_X-m_Y)\chi$. This proves all assertions for the given quasi-isometry.
\end{proof}

The sign and the normalization can also be read from the standard tree product. Give both trees unit edges, with $p$ and $q$ children per vertex. A pair of axis translations that goes up one edge in the first tree and down one in the second has $\chi=1$ in the height-compatible group. Its modular multiplier is
\[
                 e^{-(\log p-\log q)}=q/p.
\]
The same computation in the smooth model of Example~\ref{ex:models}(a) uses $\chi=t/\lambda$ and exponents $\lambda\tr D_\pm$. Their product in~\eqref{eq:modular} gives $-t(\tr D_+-\tr D_-)$, independently of the common curvature normalization.

\section{Consequences for Tukia-type theorems}\label{sec:tukia}
A Tukia-type theorem changes coordinates so that a group of uniform quasi-similarities becomes a group of similarities. In some versions amenability enters the construction of those coordinates. Theorem~\ref{thm:G} supplies a locally compact amenable acting group when the original group acts cocompactly on pairs.

Let $d$ be a continuous left-invariant $D$-homogeneous quasi-distance on $N$ inducing its topology, so $d(e^{tD}x,e^{tD}y)=e^t d(x,y)$. Any two such gauges are bi-Lipschitz equivalent by comparison on a compact unit sphere and homogeneity (use the sphere of Lemma~\ref{lem:gauge}). Replacing $D$ by $\lambda D$ replaces $d$ by $d^{1/\lambda}$. This raises both scale and distortion to the same power and keeps the action on pairs. Choose $\lambda$ as in Remark~\ref{rem:models}; the resulting gauge is equivalent to the Heintze boundary gauge. For $N\times\Q_m$ use the usual $m$-adic ultrametric and a tree of valence $m+1$ with edge length $(\log m)/\lambda$. The corresponding mixed boundary gauge is equivalent to $\max(d,d_{\Q_m})^{1/\lambda}$ \cite[Section~5.3]{DX16}.

\begin{cor}\label{cor:tukia}
Suppose $\Gamma$ acts by $L$-quasi-similarities of $(N,d)$ and cocompactly on distinct pairs, and let $Q$ be its closure in the homeomorphism group of the one-point compactification. The group $Q$ is locally compact, second countable and amenable. It acts continuously by $L$-quasi-similarities and cocompactly on pairs. Its topology is uniform convergence on compact subsets of $N$. In particular:
\begin{enumerate}
\item If $N=\R^n$ and $D$ is diagonalizable over $\mathbb C$, a bi-Lipschitz change $(\R^n,d)\to(\R^n,\widehat d_{\bar D})$ conjugates $\Gamma$ into the similarity group specified in \cite[Theorem~1.2]{DX16}, where $\bar D$ replaces each eigenvalue in Jordan coordinates by its real part and $\widehat d_{\bar D}=(d^{\mathrm{DX}}_{\bar D})^{1/\alpha}$, with $\alpha=\min\Re\operatorname{spec}(D)$ and $d^{\mathrm{DX}}_{\bar D}$ the normalized distance of \cite[Section~5.1]{DX16}. Thus the target uses the same unit-dilation convention as $d$.
\item Let $D$ be as in (1) and $m\ge2$. On $\R^n\times\Q_m$ use precisely the distance
\[
 d^{\mathrm{DX}}_{D,m}=\max(d^{\mathrm{DX}}_D,d_{\Q_m})
\]
of \cite[Sections~5.1 and~5.3]{DX16}, where $d_{\Q_m}$ is the usual $m$-adic ultrametric. If $\Gamma$ is a uniform quasi-similarity group for this distance acting cocompactly on distinct pairs, there are $\lambda>0$, an integer $s\ge2$ with $m,s$ powers of a common integer, and a bi-Lipschitz map to $(\R^n\times\Q_s,d^{\mathrm{DX}}_{\lambda\bar D,s})$ conjugating $\Gamma$ into its similarity group, as in \cite[Theorem~5.4]{DX16}.
\item For a Carnot-by-Carnot pair the conclusion of \cite[Theorem~1.2]{DFX25} holds for $\Gamma$ in the maximally symmetric normalization below. The local compactness and amenability hypothesis used for the one-dimensional ideal in \cite[Theorem~1.3]{DFX25} may also be removed in its application to the induced pair-cocompact action.
\end{enumerate}
\end{cor}
\begin{proof}
Apply Theorem~\ref{thm:G} to the indicated model. Taking limits in the original continuous gauge, as in Proposition~\ref{prop:closure}, retains the original $L$. Convergence on compact sets to a member of $Q$ bounds the scales by evaluation on a pair and bounds displacement at one point. The compactness argument of that proposition then gives convergence on the compactification. Both topologies are metrizable, so this identifies them. Since $\Gamma\le Q$, the pair action is still cocompact. For (1), the distance $d^\alpha$ is bi-Lipschitz equivalent to $d^{\mathrm{DX}}_D$. For a positive real diagonal $D$ with eigenvalues $\alpha_i$, the latter is $\max_i|x_i-y_i|^{\alpha/\alpha_i}$; real blocks give the analogous formula in complex-diagonalizable coordinates. Apply \cite[Theorem~1.2]{DX16} to $Q$ in this normalization and take the $1/\alpha$ power of both source and target distance inequalities. This gives the stated bi-Lipschitz map to $\widehat d_{\bar D}$. For (2), apply the closure argument to $d^{\mathrm{DX}}_{D,m}$ itself. Its simultaneous positive powers in both coordinates give a mixed model boundary, so Theorem~\ref{thm:G} supplies the required locally compact amenable closure. The hypotheses of \cite[Theorem~5.4]{DX16} then hold in its original distance normalization. Both conjugacies restrict to $\Gamma$. For \cite[Theorem~1.2]{DFX25}, pair cocompactness gives the radial limit point condition as in the proof of their Theorem~1.3. In that proof the amenability assumption is used for the induced pair-cocompact action, to which the same closure argument applies.
\end{proof}

\begin{rem}\label{rem:boundaryunits}
The powers in the two conclusions above keep track of the relative boundary scales. For $D=2\,\mathrm{id}$ on $\R$, the unit-dilation distance is $|x-y|^{1/2}$, whereas $d_D^{\mathrm{DX}}=|x-y|$. Their Hausdorff dimensions are $2$ and $1$, so the target in (1) must carry the $1/\alpha$ power as stated.

In the mixed case a change of normalization must be made in both coordinates. For example, on
\[
 (\R\times\Q_2,\max(|x-y|^{1/2},|u-v|_2))
\]
all translations together with $(x,u)\mapsto(4x,u/2)$ form a similarity group acting cocompactly on distinct pairs. Indeed translation puts the first point at $(0,0)$, and an integer power of this dilation puts the separation in $[1,2]$, a compact set of pair representatives. Every connected component is a snowflaked line of Hausdorff dimension $2$. In a target whose real component has the one-dimensional DX normalization, each component has Hausdorff dimension $1$, excluding a bi-Lipschitz conjugacy to that target. Part~(2) therefore specifies the complete mixed distance to which the cited theorem applies.
\end{rem}

\section{A geometric action on the original SOL-like group}\label{sec:sol}
Theorem~\ref{thm:A} supplies new Heintze factors in the smooth case. To use the original group $S$ and a metric chosen independently of $H$, we need similarity actions for prescribed homogeneous distances. This is possible for the Carnot and Carnot-by-Carnot pairs of \cite[Definition~1.1 and Section~3.4]{DFX25}. In the Carnot case the lowest layer generates the Lie algebra. In the Carnot-by-Carnot case the ideal generated by the lowest eigenspace is proper, and the quotient also carries a Carnot structure with the normalization specified there.

For such a pair $(N,D)$ let $\mu$ be the smallest eigenvalue. Use a maximally symmetric $D/\mu$-homogeneous distance $d_0$ as in \cite[Theorem~2]{HS90} and \cite[Lemma~2.3]{DFX23}. The associated inner product has orthogonal eigenspaces and is invariant under a maximal compact group $K$ of graded automorphisms. We need two consequences of this choice.

First, a second countable locally compact amenable group acting continuously by uniform quasi-similarities of $d_0$ and cocompactly on pairs can be conjugated into $\operatorname{Sim}(N,d_0)$ by a bi-Lipschitz map. For a Carnot-by-Carnot pair use \cite[Theorem~1.2]{DFX25} and the radial limit point condition above. For a Carnot group of dimension at least two apply \cite[Theorem~1.3]{DFX23} with one layer in the filtration. It yields a Carnot--Carath\'eodory target metric; maximal symmetry conjugates its similarity group into that of $d_0$ \cite[Lemma~2.3]{DFX23}. In dimension one use \cite{FM99}.

Second, a similarity of ratio $e^{\mu t}$ has the form $\sigma=L_a e^{tD}\varphi$ with $\varphi\in K$, commuting with $D$ \cite[Theorem~1.2]{KLD17}, \cite[Lemma~3.5]{DFX25}, \cite[Lemma~2.3]{DFX23}. Thus
\[
       (x,s)\longmapsto (\sigma x,s+t)=(a,t)(\varphi x,s)
\]
is an isometry for the left-invariant extension metric on $N\rtimes_D\R$ with unit vertical vector orthogonal to the chosen horizontal inner product. Evaluation and the similarity ratio give a continuous homomorphism from $\operatorname{Sim}(N,d_0)$ to these isometries.

\setcounter{mainresult}{3}
\begin{mainresult}\label{thm:D}
Let $S=(N_1\times N_2)\rtimes\R$ where $t$ acts by $(e^{tD_1},e^{-tD_2})$. Assume both factors are nontrivial, each Carnot or Carnot-by-Carnot, and $\tr D_1\ne\tr D_2$. Fix the left-invariant Riemannian metric $g_0$ with pairwise orthogonal summands, the chosen inner products on $\mathfrak n_i$, and unit vertical vector. Every compactly generated locally compact group $H$ quasi-isometric to $S$ has a compact normal subgroup $W$ such that $H/W$ is isomorphic to a closed cocompact subgroup of $\Isom(S,g_0)$. Its elements act by maps $L_s\circ(\varphi_1\times\varphi_2\times\mathrm{id})$ with $\varphi_i\in K_i$.
\end{mainresult}
\begin{proof}
Use Lemma~\ref{lem:sc}, reorder if necessary, and take the factors $X_i$ built from $\lambda D_i$ in Example~\ref{ex:models}(a). On $N_i$ their gauges are bi-Lipschitz equivalent to $d_{0,i}^{1/(\lambda\mu_i)}$. The two boundary closures $Q_i$ from Theorem~\ref{thm:projections} are amenable. By Lemma~\ref{lem:pairs} they act cocompactly on pairs. Conjugate them by bi-Lipschitz maps $F_i$ to similarities $\sigma_i(f)=F_ifF_i^{-1}$ of $d_{0,i}$. Define $\tau_i(f)$ by the similarity ratio $e^{\mu_i\tau_i(f)}$. Comparison of the two gauges and the bi-Lipschitz bounds give
\[
                \tau_i(f)=\lambda\beta_i(f)+O(1)\quad(f\in Q_i).
\]
The functions $\tau_i$ are continuous characters. Their sum on the two images of $H$ is bounded by Theorem~\ref{thm:projections}, hence zero. Write their values as $\tau,-\tau$. On the opposite-height subgroup of the product of the two extensions the resulting pair acts by
\[
        (x,y,t)\longmapsto (a_1,a_2,\tau)(\varphi_1x,\varphi_2y,t).
\]
It preserves $g_0$: the last three-coordinate map is an automorphism with orthogonal differential at the identity, and the other map is left translation.

For properness and cocompactness compare this action with the one we started with. The map $\Theta(\xi,\eta,z)\mapsto(F_1\xi,F_2\eta,\lambda z)$ from $X_1\hor X_2$ to $S$ is a quasi-isometry. The boundary changes at fixed height preserve distances up to bounded error by Lemma~\ref{lem:productdist}, and Example~\ref{ex:models}(a) supplies the remaining identification via two geometric actions of $S$. Equation~\eqref{eq:oneqi} and $\tau=\lambda b_a+O(1)$ give coarse equivariance, since a uniformly bounded vertical change along a fixed line has bounded length. Consequently the orbit is cobounded. A compact-set transporter moves a fixed base point a bounded distance, hence lies in a relatively compact subset of $H$ by comparison with the transported left translations. It is closed, so it is compact. As in Theorem~\ref{thm:modular}, these facts give a compact action kernel and a closed cocompact image in $\Isom(S,g_0)$. The compact quotient reduction lifts the same conclusion to the original $H$.
\end{proof}
The remaining condition in \cite[Theorem~1.5]{DFX25} concerns a Carnot-by-Carnot factor with a one-dimensional ideal generated by the lowest eigenspace: either the finitely generated group $\Gamma$ is amenable, or some $\Isom(S,g)$ admits a uniform lattice. The two trace cases below dispense with this condition by combining the preceding obstruction with known results.
\begin{maincor}\label{cor:E}
For a SOL-like group whose two factors are nontrivial and each Carnot or Carnot-by-Carnot, the conclusion of \cite[Theorem~1.5]{DFX25} holds without its additional assumption. In particular \cite[Conjecture~1.6]{DFX25} holds for these groups.
\end{maincor}
\begin{proof}
At unequal traces Theorem~\ref{thm:A} excludes every finitely generated group, so the conditional rigidity assertion is vacuous. At equal traces let $\Gamma$ be a finitely generated group quasi-isometric to $S$. The solvable group $S$ is amenable and unimodular, hence metrically amenable. Quasi-isometry invariance implies that $\Gamma$ is metrically amenable and therefore amenable \cite{Te08}, \cite[Lemma~19.5]{Cor18}. This verifies the amenability alternative in \cite[Theorem~1.5]{DFX25}, including when one of the indicated ideals is one dimensional. The geometric construction in \cite[Section~4.1]{DFX25} identifies each finitely generated group in the class, up to finite index and finite kernel, with a \emph{cocompact} lattice in the stated isometry group. Thus the existence of such a finitely generated group supplies a left-invariant metric whose isometry group admits a uniform lattice, which is the conclusion needed for \cite[Conjecture~1.6]{DFX25}.
\end{proof}

\section{Consequences for groups, spaces and graphs}\label{sec:applications}
\setcounter{mainresult}{1}
A closed subgroup of the isometry group of a proper geodesic space that acts cocompactly on the space is compactly generated and quasi-isometric to that space \cite[Theorem~4.C.5]{CdlH16}. Applied to a different space in the same quasi-isometry class, this gives a restriction on all its possible geometric group actions.
\begin{maincor}\label{cor:B}
If $m_X\ne m_Y$, then $X\hor Y$ is not quasi-isometric to any finitely generated group or to any unimodular or non-amenable compactly generated locally compact group. If a proper geodesic space $P$ is quasi-isometric to $X\hor Y$, every closed subgroup $G\le\Isom(P)$ acting cocompactly on $P$ is amenable, non-unimodular and contains no lattice.
\end{maincor}
\begin{proof}
A discrete group's counting measure is right invariant. Apply Theorem~\ref{thm:A}(1). For the second assertion closedness gives a proper action on $P$, and the assumed cocompactness on $P$ gives an orbit quasi-isometry. Thus Theorem~\ref{thm:A} applies to $G$. A locally compact group with a lattice is unimodular, by the invariant-measure criterion on the quotient.
\end{proof}
\begin{maincor}\label{cor:C}
\leavevmode
\begin{enumerate}
\item[(a)] Let $S=(N_+\times N_-)\rtimes\R$ be as in Example~\ref{ex:models}(a), with both factors nontrivial and unequal traces. No finitely generated group is quasi-isometric to $S$. Every compactly generated locally compact $H$ in its quasi-isometry class is amenable and non-unimodular. Modulo a compact normal subgroup it is a closed cocompact subgroup of $\Isom(S',g')$ for a SOL-like group $S'$ quasi-isometric to $S$ with a left-invariant Riemannian metric $g'$.
\item[(b)] For integers $p,q\ge2$ with $p\ne q$, every compactly generated locally compact group quasi-isometric to $\DL(p,q)$ is amenable and non-unimodular. Modulo a compact normal subgroup it acts geometrically on a product of two regular trees with opposite heights.
\item[(c)] For an integer $p\ge2$ and a real $q>1$ with $q\ne p$, every compactly generated locally compact group quasi-isometric to $\HT(q,p)$ is amenable and non-unimodular. In particular no finitely generated group is quasi-isometric to this space.
\end{enumerate}
\end{maincor}
\begin{proof}
Apply Theorem~\ref{thm:A} and Corollary~\ref{cor:B} with the rates of Example~\ref{ex:models}. In (a) the new factors are Heintze spaces. Their opposite-height hypersurface in the Riemannian product carries an induced Riemannian metric $g'$. The corresponding SOL-like group $S'$ acts simply transitively by left translations, so $g'$ is left invariant. Pairs of factor isometries preserve $g'$. Both the geometric model topology and the Riemannian topology are the hypersurface topology, so properness and cocompactness persist. A faithful proper isometric action gives a closed image topologically isomorphic to the acting group, and cocompactness makes the image cocompact in the full isometry group by the compact representative argument in Theorem~\ref{thm:modular}.
\end{proof}
Within the Carnot--Sol subclass, the isomorphism classification between such Lie groups is already \cite[Theorem~C]{Fe26b}; for Carnot and Carnot-by-Carnot factors see also \cite[Theorem~1.8]{DFX25}. Part~(a) permits an arbitrary compactly generated locally compact $H$ and is accompanied by the modular-function identity of Theorem~\ref{thm:A}.

For (a), unequal traces also force non-unimodularity of $\Isom(S,g)$ for \emph{every} left-invariant Riemannian metric $g$, by Corollary~\ref{cor:B}. This verifies the asserted non-unimodular case of \cite[Conjecture~1.6]{DFX25}. For (b) the discrete exclusion was already known \cite[Theorem~1.4]{EFW12}. For (c), \cite[Theorem~2.15 and Corollary~2.18]{BSSW15} computes the modular obstruction for a geometric group on the original treebolic model; here the exclusion applies to the whole quasi-isometry class.

\setcounter{mainresult}{5}
\begin{maincor}\label{cor:F}
If a connected locally finite quasi-transitive graph $G$ is quasi-isometric to $\DL(p,q)$ with $p\ne q$, then $\Aut(G)$ is amenable and non-unimodular. Modulo a compact normal subgroup it acts geometrically on a horospherical product of two regular trees.
\end{maincor}
\begin{proof}
In the topology of pointwise convergence vertex stabilizers are compact open. Thus $\Aut(G)$ acts properly on the metric realization, and finitely many vertex orbits give cocompactness. The orbit comparison and Corollary~\ref{cor:C}(b) give the conclusion.
\end{proof}
The Soardi--Woess criterion relates graph amenability to amenability and unimodularity of a closed transitive group \cite{SW90}. Here Corollary~\ref{cor:F} identifies both group properties in the unequal-parameter class. In particular these classes cannot supply the vertex-transitive graph with unimodular automorphism group sought in \cite[final paragraph]{Wo13}.

The hypotheses have concrete limits. At equal exponents, lattices in $\Sol$, lamplighter groups and $\mathrm{BS}(1,p)$ occur, as do the mixed examples in Example~\ref{ex:models}(d). For \cite[Question~1]{CL25} this leaves the equal-exponent case as the place to look for further discrete geometric actions with a mixed factor. If a nilpotent factor is trivial, the smooth group reduces to a Heintze group; the hyperbolic plane, for instance, is quasi-isometric to a cocompact Fuchsian group.

\appendix
\section{A boundary-measure proof of non-unimodularity}\label{sec:counting}
There is a second way to see the obstruction in Theorem~\ref{thm:A}(1). It uses the transferred quasi-action of Section~\ref{sec:product} and never needs to choose a geometric model for either boundary group. The proof turns an approximate change of height into an exactly measure-preserving action on one additional real coordinate.

\begin{lem}\label{lem:jac}
For every $(L,s)$-quasi-similarity $f$ of $\Omega_X$ and Borel set $B$ one has $\nu_X(fB)\asymp e^{m_Xs}\nu_X(B)$ with constants depending only on $X,L$.
\end{lem}
\begin{proof}
For an open set $O$ take all balls of radius at most one contained in $O$. In successive dyadic radius ranges select a maximal disjoint family among those missing all previous selections. Every ball meets a selected ball of radius at least half its radius, so enlargement of the selected balls by $2K^2$ covers $O$ by the quasi-ultrametric inequality \cite[Chapter~III]{CW71}. The selected balls are countable by positive mass and $\sigma$-finiteness. Their enlarged images have total measure at most a constant times $e^{m_Xs}\sum r_i^{m_X}$ by Ahlfors regularity. Disjointness bounds this sum by a constant times $\nu_X(O)$. Outer regularity gives the upper bound for Borel sets; use $f^{-1}$ for the lower bound.
\end{proof}

The group-theoretic input is a counting identity. Consider a window $K$ whose inverse translates see almost every point for a uniformly positive and finite amount of Haar time. In a unimodular group one can reverse the group element in this count. This bounds the measures of all iterates of the window under a commuting transformation.
\begin{lem}\label{lem:count}
Let a second countable unimodular locally compact group $H$ act by a measure-preserving Borel action on a standard $\sigma$-finite measure space $(E,\mu)$. Suppose $0<\mu(K)<\infty$ for a Borel set $K$ and
\[
    0<m_0\le N(z)=\int_H\one K(g^{-1}z)\,d\lambda_H(g)\le M_0<\infty
       \quad\text{for almost all }z .
\]
If a Borel automorphism $T$ commutes with $H$ and $\mu(TB)=r\mu(B)$ for every Borel $B$, $0<r<\infty$, then $r=1$.
\end{lem}
\begin{proof}
Tonelli's theorem and invariance give, for $B=T^nK$,
\[
 m_0\mu(B)\le\int_B N\,d\mu
 =\int_H\mu(B\cap gK)\,d\lambda_H(g)
 =\int_K\int_H\one B(gw)\,d\lambda_H(g)\,d\mu(w).
\]
The inner integrand equals $\one K(gT^{-n}w)$. Since inversion preserves Haar measure, its integral is $N(T^{-n}w)$, at most $M_0$ for almost every $w$ because $T$ preserves the measure class. Thus $m_0 r^n\le M_0$ for every integer $n$. Letting $n\to+\infty$ excludes $r>1$, and letting $n\to-\infty$ excludes $r<1$. Hence $r=1$.
\end{proof}

\begin{thm}\label{thm:counting}
If $m_X\ne m_Y$, no unimodular compactly generated locally compact group is quasi-isometric to $X\hor Y$.
\end{thm}
\begin{proof}
Assume the contrary, set $c=m_X-m_Y>0$ after ordering the factors, and pass to the second countable quotient by Lemma~\ref{lem:sc}. Use the jointly Borel data $\varphi,\bar\varphi,q_g$ of that lemma and take $o=\varphi(e)$. Write $\Omega=\Omega_X\times\Omega_Y$, $\nu=\nu_X\otimes\nu_Y$, and $u_g$ for the product of the two boundary maps of $g$; thus $u$ is a genuine Borel action. Lemma~\ref{lem:jac} applied in each coordinate yields
\begin{equation}\label{eq:jac}
     e^{cb_g-C}\nu(B)\le\nu(u_g B)\le e^{cb_g+C}\nu(B)
\end{equation}
with $C$ independent of $g,B$. Indeed the pullback product measure has the two coordinate densities as factors. It follows that $u_g$ preserves the measure class.

\emph{Choosing the logarithmic height correction.}
Write $\nu_g(B)=\nu(u_gB)$. We need derivatives with the group variable measurable as well as the boundary variable. Take $w>0$ with $\nu'=w\nu$ a probability, and finite Borel partitions increasing to the Borel $\sigma$-algebra. On an atom $P$ use $(\int_P w\,d\nu_g)/\nu'(P)$ if $\nu'(P)>0$, and zero otherwise. The numerator equals $\int (w\one P)(u_{g^{-1}}\xi)\,d\nu(\xi)$, a Borel function of $g$. For fixed $g$ these functions are conditional expectations of $d\nu_g/d\nu$ with respect to $\nu'$, by~\eqref{eq:jac}. For each fixed $g$, equation~\eqref{eq:jac} bounds the density, and all these conditional expectations on positive-measure atoms, between $e^{cb_g-C}$ and $e^{cb_g+C}$. In particular the martingale is uniformly integrable with respect to the probability $\nu'$. Martingale convergence identifies its limit with the density almost everywhere, and that limit is finite and positive. Setting the limit to one on the exceptional set gives a jointly Borel derivative $J_0(g,\cdot)$ for every fixed $g$. The exceptional set has zero $\nu'$-measure on every $g$-section, hence zero product measure by Fubini. The density bounds therefore hold for almost every $(g,\xi)$ as well.

The chain rule yields $J_0(gg',\xi)=J_0(g,u_{g'}\xi)J_0(g',\xi)$ for almost every $\xi$ for each fixed pair. To use an actual group action apply the strict cocycle theorem \cite[Theorem~B.9]{Zi84} on $(\Omega,\nu')$ to the right action $\xi\cdot g=u_{g^{-1}}\xi$ and cocycle $J_0(g^{-1},\xi)$. Returning to the left convention gives a jointly Borel $J$ on an invariant conull set, satisfying the identity there for every $g,g'$. Extend by one on its invariant complement. The strict version agrees with $J_0$ for almost every $(g,\xi)$. It is a derivative for \emph{each} $g$: for fixed $g$ choose $g'$ so that the derivative property holds both for $g'$ and $gg'$, possible by left invariance of Haar measure. Change variables by $u_{g'}$ and use the strict identity to obtain
\[
 \int_{u_{g'}B} J(g,\eta)\,d\nu(\eta)=\int_B J(gg',\xi)\,d\nu(\xi)=\nu(u_g u_{g'}B)
\]
for every $B$. In particular~\eqref{eq:jac} gives $|c^{-1}\log J(g,\xi)-b_g|\le C/c$ for almost every $(g,\xi)$.

\emph{The measure-preserving action.}
On $E=\Omega\times\R$ put
\[
 g*(\xi,t)=(u_g\xi,t+c^{-1}\log J(g,\xi)),\qquad d\mu=e^{-ct}d\nu\,dt .
\]
The strict identity gives the group law. With $\xi'=u_g\xi$ and $t'=t+c^{-1}\log J(g,\xi)$, the image measure satisfies $d\nu(\xi')=J(g,\xi)d\nu(\xi)$, whereas
\[
 e^{-ct'}=e^{-ct}J(g,\xi)^{-1},\qquad dt'=dt
\]
for each fixed $\xi$. The factors cancel, so $\mu(g*B)=\mu(B)$ for every Borel set $B$. Vertical translation $T(\xi,t)=(\xi,t+1)$ commutes with $H$ and scales $\mu$ by $e^{-c}$.

\emph{A finite-measure window.}
By~\eqref{eq:oneqi} the corrected action satisfies
\[
                  d_Z(\Theta(g*z),q_g(\Theta z))\le A
\]
for almost every $(g,z)$, with one constant $A$. Let $U$ be a symmetric compact generating neighbourhood, choose $R>A+1+\sup_{a\in U}d_Z(\varphi(a),o)$, and set $K=\Theta^{-1}(\overline B(o,R))$. Lemma~\ref{lem:productdist} bounds the height on $K$ and puts both boundary coordinates in fixed compact balls. The window is closed and contains a neighbourhood of a preimage of $o$. Hence $0<\mu(K)<\infty$.

For almost every $z$ the comparison holds with $g^{-1}$ for almost every $g$, by Fubini and inversion invariance. Set $v=\bar\varphi(\Theta z)$. For $g\in vU$, $g^{-1}v\in U$, so $q_{g^{-1}}(\Theta z)=\varphi(g^{-1}v)$ lies within $R-A$ of $o$. Thus the count in Lemma~\ref{lem:count} is at least $\lambda_H(U)$. Conversely, if $g^{-1}*z\in K$ then, apart from the exceptional null set, $d_Z(\varphi(g^{-1}v),o)\le R+A$. The lower quasi-isometry bound puts $g^{-1}v$ in a fixed symmetric word ball $B_H(e,R')$, hence $g\in vB_H(e,R')$. This bounds the count above by the finite Haar measure of that ball. Lemma~\ref{lem:count} now gives $e^{-c}=1$, contrary to $c>0$.
\end{proof}
For a general locally compact group inversion contributes a modular weight to the count. In the smooth model the vertical translation of the extended space is right multiplication and has precisely that weight, as predicted by~\eqref{eq:modular}.

\section*{Data availability}
No datasets were generated or analysed during this study.

\section*{Acknowledgements}
This work was supported by the National Natural Science Foundation of China (NSFC), grant no.~12131012.

\section*{Declaration of generative AI and AI-assisted technologies in the manuscript preparation process}
During preparation of the manuscript, the authors used OpenAI's ChatGPT (Codex) to improve the English presentation, and expand intermediate estimates. The tool also suggested strengthening the non-unimodularity conclusion to the modular-function identity in Theorem~\ref{thm:A} and assisted in developing and drafting its proof in Section~\ref{sec:mod}. The tool also assisted in assessing a mathematical review report and in revising the affected statements and calculations. The authors take responsibility for the mathematical arguments, references, and text.

\end{document}